\documentclass[12pt]{iopart}

\usepackage{multirow}
\usepackage{longtable,tabularx}
\usepackage{graphicx} 

\expandafter\let\csname equation*\endcsname\relax
\expandafter\let\csname endequation*\endcsname\relax
\usepackage{amsmath,bm}

\usepackage{amsfonts}
\usepackage{enumitem}
\usepackage{xcolor}
\usepackage{algorithm}
\usepackage[noend]{algpseudocode}
\algnewcommand\algorithmicinput{\textbf{Input:}}
\algnewcommand\Input{\item[\algorithmicinput]}
\algnewcommand\algorithmicoutput{\textbf{Output:}}
\algnewcommand\Output{\item[\algorithmicoutput]}

\usepackage{nomencl}
\makenomenclature

\usepackage{longtable,tabularx}
\usepackage{longtable,tabularx}
\begin{document}

% \title[ ]{The albatross optimized flight as an extremum seeking phenomenon: autonomous and stable controller for dynamic soaring}

\title[]{A novel hypothesis for how albatrosses optimize their flight physics in real-time: an extremum seeking model and control for dynamic soaring}
\author{Sameer Pokhrel}
\address{PhD Student\\Department of Aerospace Engineering and Engineering Mechanics,\\University of Cincinnati, OH, USA\\
Email: pokhresr@mail.uc.edu}

\author{Sameh A. Eisa}
\address{Assistant Professor \\Department of Aerospace Engineering and Engineering Mechanics,\\University of Cincinnati, OH, USA\\
Email: eisash@ucmail.uc.edu}

% \address{IOP Publishing, Temple Circus, Temple Way, Bristol BS1 6HG, UK}
% \ead{submissions@iop.org}
% \vspace{10pt}
% \begin{indented}
% \item[]August 2017
% \end{indented}

\begin{abstract}
The albatross optimized flight maneuver -- known as dynamic soaring -- is nothing but a wonder of biology, physics, and engineering. By utilizing dynamic soaring, this fascinating bird can travel in the desired flight direction almost for free by harvesting energy from the wind. This phenomenon has been observed for centuries as evidenced by the writings of Leonardo da Vinci and Lord Rayleigh. Moreover, dynamic soaring biological inspiration has triggered a momentous interest among many communities of science and engineering, particularly aeronautical, control, and robotic engineering communities. That is, if dynamic soaring is mimicked, we will have arrived at a new class of unmanned aerial vehicles that are very energy-efficient during part (or the full) duration of their flight. Studying, modeling, and simulating dynamic soaring have been conducted in literature by mostly configuring dynamic soaring as an optimal control problem. Said configuration requires accurate dynamic system modeling of the albatross/mimicking-object, accurate wind profile models, and a defined mathematical formula of an objective function that aims at conserving energy and minimizing its dissipation; the solution then of such optimal control problem is the dynamic soaring trajectory taken -- or to be taken -- by the bird/mimicking-object. Furthermore, the decades-long optimal control configuration of the dynamic soaring problem resulted in non-real-time algorithms and control solutions, which may not be aligned well with the biological phenomenon itself; experimental observations of albatrosses indicate their ability to conduct dynamic soaring in real-time. Indeed, a functioning modeling and control framework for dynamic soaring that allows for a meaningful bio-mimicry of the albatross needs to be autonomous, real-time, stable, and capable of tolerating the absence of mathematical expressions of the wind profiles and the objective function -- hypothetically similar to what the bird does. The qualifications of such modeling and control framework are the very same characteristics of the so-called extremum seeking systems. In this paper, we show that extremum seeking systems existing in control literature for decades are a natural characterization of the dynamic soaring problem. We propose an extremum seeking modeling and control framework for the dynamic soaring problem hypothesizing that the introduced framework captures more features of the biological phenomenon itself and allows for possible bio-mimicking of it. We provide and discuss the problem setup, design, and stability of the introduced framework. Our results, supported by simulations and comparison with optimal control methods of the literature, provide a proof of concept that the dynamic soaring phenomenon can be a natural expression of extremum seeking. Hence, dynamic soaring has the potential to be performed autonomously and in real-time with stability guarantees.
\end{abstract}

% \nomenclature{\(DS\)}{Dynamic Soaring}
% \nomenclature{\(ESC\)}{Extremum Seeking Control}
% \nomenclature{\(UAVs\)}{Unmanned Aerial Vehicles}
% \nomenclature{\(ESC\)}{Extremum Seeking Control}
% \nomenclature{\(x\)}{Position vector along East direction}
% \nomenclature{\(y\)}{Position vector along North direction}
% \nomenclature{\(z\)}{Altitude}
% \nomenclature{\(V\)}{Airspeed}
% \nomenclature{\( \gamma\)}{Air relative flight path angle}
% \nomenclature{\(\psi \)}{Air relative heading angle}
% \nomenclature{\( \phi\)}{Bank angle}
% \nomenclature{\(C_L\)}{Coefficient of Lift}
% \nomenclature{\(C_D\)}{Coefficient of Drag}
% \nomenclature{\(L\)}{Aerodynamic Lift}
% \nomenclature{\(D\)}{Aerodynamic Drag}
% \nomenclature{\(\rho\)}{Air density}
% \nomenclature{\(S\)}{Wing area}
% \nomenclature{\(m\)}{Mass of UAV/albatross}
% \nomenclature{\(n\)}{Load factor}
% \nomenclature{\(SISO\)}{Single input single output}
% \nomenclature{\(K\)}{Induced drag coefficient}
% \nomenclature{\(W_0\)}{Free stream wind speed}
% \nomenclature{\(\delta\)}{Shear layer thickness}
% \nomenclature{\(PE\)}{Potential energy}
% \nomenclature{\(KE\)}{Kinetic energy}
% \nomenclature{\(TE\)}{Total energy}
% \nomenclature{\(DC\)}{Direct collocation}
% \nomenclature{\(W\)}{Wind velocity}
% \nomenclature{\(\dot{W}\)}{Wind shear gradient}
% \nomenclature{\(T_{cyc}\)}{Time period for a single cycle}
% \nomenclature{\(TE\)}{Total Energy}
% \nomenclature{\(g\)}{Acceleration due to gravity}
% \printnomenclature
\MakeUppercase{List of important symbols and abbreviations as they appear on the text}
{\renewcommand\arraystretch{1}
\begin{longtable}{@{}l @{ \quad} l@{}}
$ESC$&Extremum Seeking Control\\
$GPOPS2$&General purpose optimal control software\\
$W_0$&Free stream wind speed for logistic wind model\\
$W$&Wind velocity\\
$z$&Altitude\\
$z_m$& Altitude corresponding to the middle of shear layer for logistic wind \\
$\delta$&Shear layer thickness for logistic wind model\\
$z_{ref}$ & {Reference altitude for logarithmic wind model}\\
{$V_{w_{ref}}$ }& {Value of wind shear strength at reference altitude for logarithmic wind model }\\
{$z_0$ }& {Surface correction factor for logarithmic wind model}\\
%10 done
$x$&Position vector along East direction\\
$y$&Position vector along North direction\\
$V$&Airspeed\\
$ \gamma$&Air relative flight path angle\\
$\psi $&Air relative heading angle\\
$ \phi$&Bank angle\\
$L$&Aerodynamic Lift\\
$D$&Aerodynamic Drag\\
$\dot{W}$&Wind shear gradient\\
$m$&Mass of UAV/albatross\\
%20 done
$g$&Acceleration due to gravity\\
$C_L$&Coefficient of Lift\\
$C_D$&Coefficient of Drag\\
$\rho$&Air density\\
$S$&Wing area\\
$C_{D0}$ & Zero-lift drag coefficient\\
$K$&Induced drag coefficient\\
$J$& Performance Index\\
$t_0$& Initial time for dynamic soaring cycle\\
$t_f$& Final time for dynamic soaring cycle\\
%30 done
$T$& Duration for dynamic soaring cycle\\
$(\cdot)_{min}$ & Minimum value\\
$(\cdot)_{max}$ & Maximum value \\
$n$&Load factor\\
$a$ & Amplitude of modulation signal\\
$b$ & Amplitude of demodulation signal\\
$\omega$ & Frequency of input signals\\
$\phi_{phase}$& Phase lag between modulation and demodulation signals\\
$L[\cdot]$& Laplace transform\\
% 41 completed
$\eta_{dist}$& Added disturbance/noise\\
$SISO$&Single-input single-output\\
$TE$ &Total Energy\\
$e$& Specific total energy\\
$PE$&Potential energy\\
$KE$&Kinetic energy\\
$\eta,\xi$ & Intermediate states for classic extremum seeking control\\
$F_i,F_0$& Input and output dynamics for system for augmented extremum seeking\\
$C_i,C_0$& Compensators used in design of augmented extremum seeking structure\\
% {$int()$}& { Interior}\\
% {$conv()$}&{ Convex hull}\\
% {$U$ }&{Set of admissible controls}
\end{longtable}}

\section{Introduction}
\begin{figure}[ht]                          % [htb]% order of placement preference: here, top, bottom
    \centering
	\includegraphics[width=\textwidth]{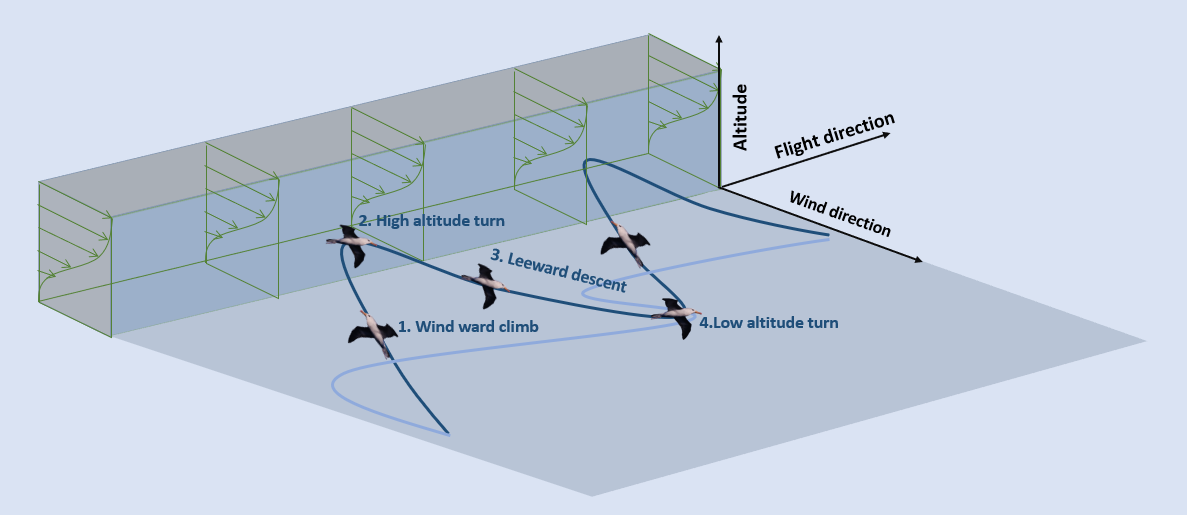}
	\caption{Dynamic soaring maneuvers conducted by albatrosses or other soaring birds are represented by a black trajectory in the presence of wind spatial distribution depicted in green. Wind speed changes considerably with height causing what is known as wind shear/gradient which is a prerequisite for dynamic soaring. The bird goes into the headwind to acquire lift. The acquired lift from the wind then allows the bird to gain height up to a point. The bird then turns and glides/descends trading the gained potential energy with kinetic energy causing gains in velocity and momentum. At a low altitude, energized by the gained velocity/momentum, the bird then turns into the headwind and repeats the dynamic soaring cycle again.}
	\label{fig:dynamicsoaring}
\end{figure}

Dynamic soaring is a fascinating flight maneuver that enables long-duration flights by harvesting atmospheric energy \cite{gao2017dubins,mir2019soaring}. The energy needed to perform such a long duration of flights is gained from the wind in the proximity of the surface. In regions above the sea surface or mountainous areas, the speed of the wind changes considerably with the altitude to yield what is known as ``wind shear/gradient" \cite{sachs2013experimental,mir2018review}. By flying across the wind gradient region periodically, energy is harvested from the spatial wind speed distribution. This has been observed among what is known as soaring birds such as albatrosses and eagles, among others. As a matter of fact, the dynamic soaring maneuver enables soaring birds to travel large distances almost without flapping their wings as energy spending is minimized substantially \cite{richardson2011albatrosses}. The dynamic soaring maneuver can be considered a cycle of four characteristic flight phases \cite{mir2018review}: (i) windward climb, (ii) high altitude turn, (iii) leeward descent, and (iv) low altitude turn -- see figure \ref{fig:dynamicsoaring}. In order to conduct this maneuver, the bird goes into the headwind and gains height -- trading off kinetic energy with potential energy. The lift acquired from the wind allows the bird to sustain gains in the height up to a point. Then, the bird takes a steep turn and dives down (glide) with a tailwind. It continues to descend -- trading off potential energy with kinetic energy, which then transfers into gains in the velocity. At a low altitude point, the bird, energized by the gained velocity (momentum), takes the low altitude turn and starts a new dynamic soaring cycle. Ideally, the dynamic soaring cycle is an energy-neutral (near-neutral) maneuver. It is important to emphasize that, the dynamic soaring phenomenon has been verified and validated experimentally \cite{sachs2013experimental,yonehara2016flight}.

The interest in the albatross flight secret goes back to Leonardo da Vinci and Lord Rayleigh \cite{richardson2019leonardo} hundreds of years ago. This interest is still present today among many communities of science and engineering \cite{denny2008dynamic,wilson1975sweeping,rayleigh1883soaring,sachs2013experimental}. From physical and engineering points of view, dynamic soaring is an extremely rare -- if not entirely unique flight dynamic system. This is due to the fact that in flight dynamic systems, a price has to be paid for the drag (the air resistance and friction). Consequently, flight dynamic systems are not conservative systems as they dissipate energy to overcome the drag. Dynamic soaring system on the other hand, performs as a conservative-like system, i.e., the sum of its kinetic and potential energy remains near-constant without dissipation of energy, in contrast to almost all other flight dynamic systems. Not surprising then, the potential of dynamic soaring bio-inspiration has captured exceptional interest in the aeronautical, robotic, and control engineering communities. That is, if one will be able to mimic soaring birds, then substantial advancements can be projected on unmanned aerial vehicles and bio-inspired robotics. In order to do such a complex bio-mimicry, researchers tried to study the modeling and control methods, designs, and techniques that will allow the simulation, mimicking, and eventually, the application of dynamic soaring in real life, especially to unmanned systems and drones \cite{zhao2004optimal,deittert2009engineless,sukumar2013sailplanes, bird2014closing}. However, as documented in the review work \cite{mir2018review},  most of the modeling, simulation, and control literature on dynamic soaring are configured as an optimal control problem, in part or in full. That is, the objective of the dynamic soaring problem is to solve for, or find the, optimal trajectory, through which the flight is energy neutral, or near-neutral \cite{sachs2005minimum,mir2018optimal,gao2017dubins,MITbousquet2017optimal}. It is important to emphasize that finding such optimal trajectory -- the solution of the dynamic soaring optimal control problem -- typically is not achievable in real-time throughout the literature \cite{mir2018review,mir2018optimal}. Moreover, said optimal control solutions of dynamic soaring are based on algorithms and methods that are very dependent on the mathematical formulated model/expression of the wind shear and the objective function that is aiming at minimizing energy spending or maximizing energy gain. Also, dynamic soaring problem setup in all the above-mentioned contributions in literature usually necessitates heavy constraints and bounds, which have to be imposed \cite{mir2018review}. Recently, a dynamic soaring simulation
system with distributed pressure sensors, that contains an online estimation and control
stage utilizing an offline training framework, has been developed \cite{wangTrainingDS2022bio}. In that work, wind information can be estimated online using the trained model and then is implemented in the simulation of dynamic soaring. 
Nevertheless, said recent study is not operable in real-time as well. 
% It is worth noting here that from a control theory and systems technical point of view, in literature, there is no access to a real-time control design for mimicking dynamic soaring with stability guarantees} \cite{mir2021stability}.

Motivation: Our motivation is to introduce and identify an alternative control theory and methods that capture better the nature of the dynamic soaring problem. This will enable novel modeling, simulation, and control tools for the dynamic soaring phenomenon. It is quite clear that the control methods of the dynamic soaring problem in literature lead to systems that are one or more of the following: (i) non-real-time, (ii) computationally complex and expensive, (iii) unstable, and (iv) very model-dependent (the objective function and wind profile models have to be known a priori). By looking at the literature on autonomous control systems, one class of systems stands out: Extremum Seeking Control (ESC) systems \cite{TanHistory2010}. These systems seem to be very descriptive and natural to the dynamic soaring problem. Furthermore, they do not possess the above-mentioned flaws, points (i)-(iv). These extremum seeking systems operate by steering a dynamical system -- such as the albatross or a mimicking-object -- to the extremum (maximum or minimum) of an objective function. This objective function can be, for example, maximum energy gain or minimum energy spending, among many other objectives. Moreover, extremum seeking systems admit objective functions without requiring access to their formulated mathematical expression as long as measurements of such objective functions can be obtained. Furthermore, and perhaps the most powerful trait of extremum seeking systems, they can tolerate less accurate modeling, or no model at all, of the dynamical system itself, as long as the objective function can be measured accurately. That is, in summary, extremum seeking systems are: (i) model-free, (ii) autonomous, (iii) real-time, and (iv) stable \cite{ariyur2003real,scheinker2017model}. Of note, the particular extremum seeking structures investigated in this paper in later sections have been thoroughly and rigorously characterized in \cite{KRSTICMain,krstic2000performance,ESC2_2}. Extremum seeking systems, as described above, seem to be strongly fitting to bio-mimicking modeling in general and the dynamic soaring problem in particular as the bird conducts dynamic soaring certainly without comprehending the mathematical closed forms of the wind profile or the objective function, and it definitely conducts the dynamic soaring flight in real-time \cite{sachs2013experimental} while sensing (measuring) its surrounding. For instance, the nostrils of the albatross act as an airspeed sensor \cite{pennycuick2008information_nostrils,brooke2002gusts_nostrils}. 

Contribution: In this paper, we hypothesize that extremum seeking systems found in the literature for decades \cite{TanHistory2010,ariyur2003real} are able to describe, model, capture, and replicate/mimic the dynamic soaring phenomenon that is conducted by the albatross and other soaring birds. We construct our proof of concept for said hypothesis by providing: (i) a novel characterization of the energy neutral (near neutral) dynamic soaring maneuver/cycle as an autonomous, real-time, extremum seeking system model; (ii) successful implementation of two extremum seeking control structures incorporating the dynamic soaring problem; (iii) stability analysis of extremum seeking systems characterizing dynamic soaring even under some disturbance/noise in the objective function measurement; we provide a rigorous mathematical theorem for stability as well; and (iv) a comparison between the results of the new implemented systems against the solution of the powerful optimal/numerical optimizer: GPOPS2 \cite{gpops2} representing the decades-long optimal control literature of dynamic soaring. In addition to the contributions directly related to the modeling and control of the dynamic soaring problem, it is important to note that, this work, to the best of our knowledge, is the first to bring extremum seeking systems and their respective theory and methods to the field of bird flights and their bio-inspiration. Extremum seeking systems are appealing for characterizations in biological systems and applications to bio-mimicry due to their ability to describe naturally the steering of systems to optimal state (extremum of an objective function) as done, arguably, by many animals, birds, and insects based on sensory (measurement) of different signals, for example, sensing light (phototaxis), chemical concentration (chemotaxis),
temperature (thermotaxis), etc. \cite{dusenbery1992sensory_taxis,newman1982infrared_snake_thermotaxis}. We believe this paper can encourage, and give insights, for more studies and lines of research on discovering and revealing extremum seeking systems in nature with novel modeling and control applications.

Organization of the paper: In section \ref{problem_setup}, we provide some background on the decades-long literature of configuring the dynamic soaring problem as an optimal control problem, in part or in full. We also provide in the same section, the description, problem setup, the models of the system dynamics, and the wind shear profiles used in literature and in some simulations in this paper. In section \ref{ESC}, we provide a brief background on extremum seeking systems to especially help readers who are not familiar with this kind of systems and controls. In the same section, we provide a couple of examples and the main control designs/structures related to this paper. In section \ref{ESC_DS}, we introduce our proposed extremum seeking model for dynamic soaring and two control structures to solve the problem. In the same section, we provide our simulation results with a comparison against the powerful optimal control numerical solver, GPOPS2 \cite{gpops2}. We also provide a thorough discussion of our results and stability analysis at the end of that section. Section \ref{conclusion} concludes the paper and provides a future prospective.
\section{Dynamic soaring in literature: an optimal control problem configuration}\label{problem_setup}
\begin{figure}[ht]
    \centering
    \includegraphics[width=0.75\textwidth]{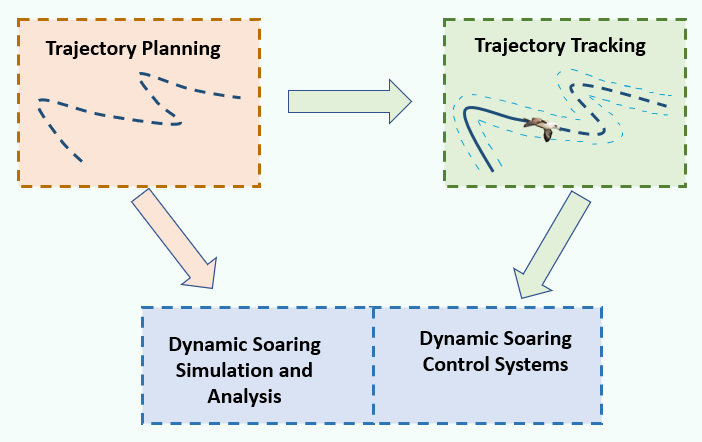}
    \caption{Two lines of research are commonly found in dynamic soaring literature. The first line, depicted in orange, consists of studies related to trajectory planning, generation, simulation, and the analysis of such trajectories. The dashed blue curve in the orange block is the planned trajectory based on the numerical solution of the optimal control configuration; the bird/mimicking-object is projected to take (have taken) this trajectory which is computed in anytime between 8 to 1000 seconds (typically not real-time). The second line of research, depicted in green, involves trajectory planning as well as developing controllers for trajectory tracking; hence, creating a dynamic soaring control system. Studies in this line of research still have to obtain the planned trajectory from studies in the orange part. After that, a control architecture/design is implemented to have the system, e.g., a mimicking object, tracking the planned trajectory within a tolerable error. In the green block, the already tracked part of the trajectory is depicted in a solid blue curve, and the rest of the ``to be tracked" part of the trajectory is depicted in a dashed blue curve. The error margin for the tracking controller is depicted in lighter dashed curves about the trajectory of interest. }
    \label{fig:traj_plan}
\end{figure}

Recall from the introduction, dynamic soaring in literature \cite{mir2018review} is configured typically as an optimal control problem, in part or in full. That is, researching dynamic soaring in literature can be broadly categorized into two lines of research. The first line consists of contributions that are focused on trajectory planning, generation, simulation, and the analysis of such trajectories. For example, finding trajectories for various modes of dynamic soaring has been studied in \cite{zhao2004optimal}, finding trajectories corresponding to minimum wind strength, maximum speed, peak altitude, etc., have been studied in \cite{sachs2005minimum,pennycuick1982flight}, and finding robust optimal trajectory is studied in \cite{flanzer2012robust}. Furthermore, analysis of the optimal trajectory from controllability and stability perspectives is also done in \cite{mir2018controllability, mir2021stability}. Most of the studies in this line of research utilize iterative numerical optimization techniques such as advanced launcher trajectory optimization software (ALTOS), graphical environment for simulation and optimization (GESOP), nonlinear programming solver (NPSOL), imperial college London optimal control software (ICLOCS), inverse dynamics in the virtual domain (IDVD), and general purpose optimal control software (GPOPS) \cite{sachs2003optimization_software,sachs2005minimum,sachs2001shear_software,zhao2004optimal,gpops2}, with computational time ranging from 8 seconds to 1000 seconds depending on the machine, complexity of the problem, desired accuracy, and the utilized algorithm \cite{mir2018optimal}. In short, these numerical methods are utilized to solve the optimal control configuration of dynamic soaring by determining (typically not in real-time) the optimal dynamic soaring trajectory per the configuration. The second line of research augments the first line of research by adding a control design component to the research. In other words, the second line of research is consistent with contributions aiming at developing controllers for trajectory tracking of the optimal dynamic soaring trajectory, thus creating a dynamic soaring control system. For tracking of optimal trajectories, a feed-forward-feedback control architecture is utilized in \cite{bird2014closing}, a custom-built nonlinear control technique is utilized in \cite{gao2017dubins} and a bio-inspired controller using fuzzy rules is developed in \cite{barate2006designfuzzy}. Similarly, in \cite{li2020parameterized}, authors use cubic splines and skewed/flattened sinusoids for parameterized trajectory planning for dynamic soaring. For tracking the planned trajectory, they use a combination of feedback and a feed-forward controller. The two lines of research found in the literature for dynamic soaring have been visually depicted in figure \ref{fig:traj_plan}. The orange part represents the first line of research, whereas the green part represents the second line of research.

Next, we provide the dynamic soaring optimal control problem configuration as usually done in literature \cite{mir2018review}. We believe this will be useful and important for this paper and the reader in multiple aspects: (i) configuring dynamic soaring as an optimal control problem has been the most common and nearly a singular method pursued in the literature of the problem, so it will be important for interested readers to have a more detailed, yet brief, understanding of such configuration; (ii) the main contribution of this paper is to provide an alternative to this decades-long optimal control configuration, so it is important to show contrasts and commonalities between our extremum seeking novel approach, provided in section \ref{ESC_DS},  and the optimal control configuration, while involving/citing directly elements of this section in the discussion; and (iii) we use the optimal control numerical solver GPOPS2 \cite{gpops2} as a representative for optimal control solutions when we compare simulation results between extremum seeking systems and the literature in section \ref{simulationSection}, so the reader, if interested, should be able to understand how this solver computes the optimal control solution.

\subsection{Wind shear model} 
\begin{figure}[ht]
    \centering
    \includegraphics[width=0.75\textwidth]{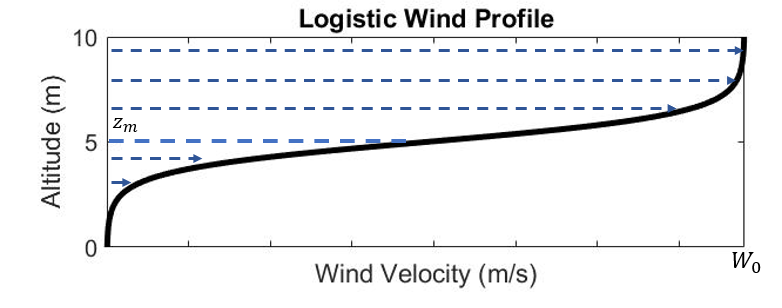}
    \caption{A logistic wind profile with a wind inhomogeneity $W_0 =7.8\: m/s$, developing over a length $\delta =2/3\:m$ with the altitude corresponding to the middle of the shear layer $z_m =5\:m$. The blue arrows show the intensity of wind at a given altitude and the black curve shows the nature of wind velocity with altitude.  This profile provides a way to capture the features of a wide range of wind fields using the three parameters $W_0$, $\delta$, and $z_m$.}
    \label{fig:sigmoidWind}
\end{figure}

\begin{figure}[ht]
    \centering
    \includegraphics[width=0.75\textwidth]{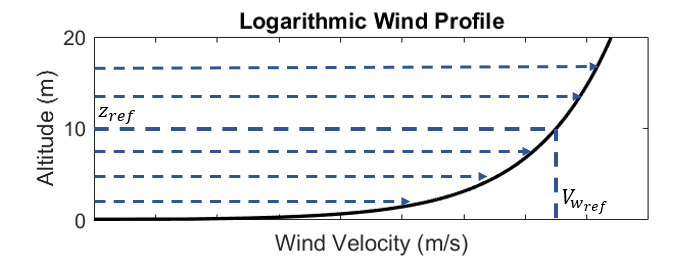}
    \caption{A logarithmic wind profile with the wind shear strength of $V_{W_{ref}}=15\: m/s$, at the reference altitude of $z_{ref}=10\: m$. The blue arrows show the intensity of wind at a given altitude and the black curve shows the nature of wind velocity with altitude. This model is commonly used for meteorological study and has been used as a wind model to study dynamic soaring in many references including \cite{sachs2005minimum,mir2018review}.}
    \label{fig:LogWindShear}
\end{figure}

Wind shear/gradient takes place on thin layers between two regions in the atmosphere if the airflow vector is different, i.e., the wind speed differs significantly between these two regions. Wind shear is a prerequisite for energy extraction in dynamic soaring \cite{mir2018review,MITbousquet2017optimal,mir2018controllability,goto2022extinct_bird_sigmoid}. Therefore, it is important to have a proper model to describe the wind dynamics for successful implementation when using dynamic optimization methods and numerical optimizers applied to dynamic soaring \cite{mir2018review}. We use two different models to describe the wind shear above sea level. The first model is a logistic model similar to
\cite{MITbousquet2017optimal,bird2014closing,wangTrainingDS2022bio}:
\begin{equation}\label{eqn:wind}
    % W(z)=\frac{W_0}{1+e^{-z/\delta}}.
    W(z)=\frac{W_0}{1+e^{-(z-z_m)/\delta}},
\end{equation}
where $z$ represents altitude. This wind model (graphically depicted in figure \ref{fig:sigmoidWind}) is parameterized by free stream wind speed $W_0$, shear layer thickness $\delta$, and the altitude corresponding to the middle of the shear layer $z_m$.
% It captures the main features of separated winds over ocean waves as well as
It captures the features of a wide range of wind fields with a typical wind inhomogeneity, $W_0$, developing over a length $\delta$ \cite{MITbousquet2017optimal}.
The second model is a logarithmic wind model, similar to \cite{sachs2005minimum,mir2018review}, given by
\begin{equation}\label{eqn:logwind}
    W(z)=V_{W_{ref}}\frac{ln(z/z_0)}{ln(z_{ref}/z_0)},
\end{equation}
where $z$ is the altitude, $V_{W_{ref}}$ is the value of wind shear strength at the altitude of $z_{ref}$, and $z_0$ is the surface correction factor that reflects surface properties like irregularities, roughness, and drag. Figure \ref{fig:LogWindShear} graphically shows the logarithmic wind profile.

\subsection{Flight dynamics model}
\begin{figure}[ht]
    \centering
    \includegraphics[width=0.75\textwidth]{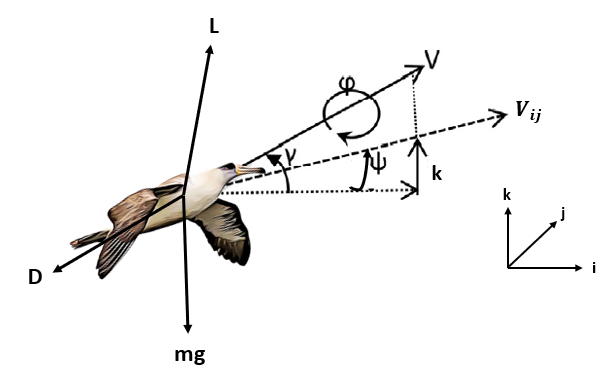}
    \caption{Aerodynamic forces lift (L) and drag (D) acting on the bird of weight $mg$, where $m$ is the mass of the bird and $g$ is the acceleration due to gravity. The angle $\psi$ is the angle between $i$ and projection of velocity $V$ in $ij$-plane, and $\gamma$ is the angle between $V$ and $ij$-plane with nose up as positive, and $\phi$ is the roll angle. Our formulation follows East, North and Up frame of reference, i.e. $(i,j,k)=(e_{East}, e_{North}, e_{Up})$}
    \label{fig:Angles}
\end{figure}

Here we adopt the three degrees of freedom, point-mass model, with no thrust component, i.e., no engine or source of power on board is assumed in the system. This is usually how an albatross or a mimicking unmanned system conducting dynamic soaring is modeled in literature \cite{mir2018review}. It is worth noting that some researchers have used six degrees of freedom models \cite{6DOFakhtar2012positioning}, however, as assessed in the review work \cite{mir2018review}, increased accuracy using six degrees of freedom models would not substitute the increased computational complexity and expense when compared with the three degrees of freedom models. Furthermore, since soaring birds such as albatrosses barely flap their wings while conducting dynamic soaring \cite{mir2018optimal,catry2004sustained_noflap,warham1996behaviour_noflap,jouventin1990satellite_noflap}, then point mass modeling of albatross/mimicking-object is considered acceptable. The model used in this work is represented by the following system of differential equations \cite{MITbousquet2017optimal}:
\begin{equation}\label{eqn:dynamics}
    \begin{split}
        \dot{x}&=V \cos\gamma \cos\psi,\\
        \dot{y}&=V \cos\gamma \sin\psi-W,\\
        \dot{z}&= V \sin\gamma,\\
        m \dot{V}&= -D -mg\sin\gamma+m\dot{W} \cos\gamma \sin\psi,\\
        mV\dot{\gamma} &= L \cos \phi - mg \cos \gamma - m \dot{W} \sin \gamma \sin \psi, \\
        mV \dot{\psi }\cos \gamma  &= L \sin \phi + m \dot{W} \cos \psi, \\
    \end{split}
\end{equation}
where $m$ represents mass, $g$ represents acceleration due to gravity and $x,y,z$ represent the position of the albatross/mimicking-object. 
The above model follows the East, North and Up frame of reference, i.e. $(i,j,k)=(e_{East}, e_{North}, e_{Up})$ as shown in figure \ref{fig:Angles}. Here, $\psi$ is the angle between $i$ and projection of velocity $V$ in $ij$-plane, and $\gamma$ is the angle between $V$ and $ij$-plane with nose up as positive and $\phi$ is the roll angle. We assume the wind is blowing from North to South (recall figure \ref{fig:dynamicsoaring}) all the time and only a horizontal wind component exists. In this model, the speed and angles of the albatross/mimicking-object are modeled in wind relative reference; and position $(x,y,z)$ is modeled in Earth fixed-frame. We calculate the lift and drag forces using the following equations:
\begin{equation}\label{DS_general_system}
    \begin{split}
        L&= \frac{1}{2}\rho V^2 S C_L,\\
        D&= \frac{1}{2}\rho V^2 S C_D,
    \end{split}
\end{equation}
where $S$ is the wing area of the bird, $\rho$ is the air density, and $C_D$ is the parabolic drag coefficient, which is given by $C_D = C_{D0} + KC_L^2$. The terms $C_{D0}$ and $K$ represent zero-lift drag coefficient and induced drag coefficient, respectively. Table \ref{tab:params} provides the albatross flight dynamic characteristics and table \ref{tab:Windparams} provides parameters used to model the wind shear.
 \begin{table}[h]
\centering
\resizebox{0.4\textwidth}{!}{
    \begin{tabular}{|>{\arraybackslash}p{2cm}|>{\arraybackslash}p{3cm}|}
		\hline 
		 Parameter & Value\\
		 \hline
		{m} &{ 8.5 kg}\\
	    {S} & {0.65 $m^2$}\\
	    {$C_{D0}$}&{ 0.033}\\
	    {K} &{ 0.019}\\
		{$\rho$} &{ 1.225 $kg/m^3$}\\
		{g }&{ 9.8 $m/s^2$}\\
		\hline
	\end{tabular}
	}
	\caption{Parametric characteristics of the albatross and the environment used.}
	\label{tab:params}
\end{table}

 \begin{table}[h]
\centering
\resizebox{0.4\textwidth}{!}{
    \begin{tabular}{|>{\arraybackslash}p{2cm}|>{\arraybackslash}p{3cm}|}
		\hline 
		 Parameter & Value\\
		 \hline
		{$W_0$} &{ 7.8 m/s }\\
		{$\delta$} &{ 2/3 m}\\
		{$z_m$} &{ 5 m}\\
		{$V_{W_{ref}}$} & {15 m/s}\\
		{$z_{ref}$} & {10 m}\\
		{$z_0$} &{ 0.03 m}\\
		\hline
	\end{tabular}
	}
	\caption{Wind model parameters used in this study.}
	\label{tab:Windparams}
\end{table}

\subsection{Problem formulation}
Now we couple the flight dynamic system model (\ref{eqn:dynamics}) subject to the wind conditions identified by the wind shear models (\ref{eqn:wind}) and (\ref{eqn:logwind})  with a mathematically defined performance index (objective function); hence, we furnish the dynamic soaring nonlinear optimal control problem formulation. 

We define the state vector $\bm{x}(t)$ and the control vector $\bm{u}(t)$ as 
\begin{equation}\label{eqn:probForm}
\begin{split}
    \bm{x}(t)&=[x,y,z,V,\gamma,\psi],\\
    \bm{u}(t)&=[C_L,\phi].
\end{split}
\end{equation}
where the states represent the velocity, position, and orientation of the bird/mimicking-object, dynamics of which is given by (\ref{eqn:dynamics}). The control inputs are the roll angle $\phi$ and coefficient of lift $C_L$, which relates the lift generated by the bird to air density, velocity, and the surface area of the wing. The flight dynamic model in (\ref{eqn:dynamics}) with the wind profile in (\ref{eqn:wind}) and (\ref{eqn:logwind}) can be configured as a nonlinear system of the form 
\begin{equation}
    \bm{\dot{x}}(t)=\bm{f}(\bm{x}(t),\bm{u}(t)).
\end{equation}

The performance index to be maximized (or minimized) - the objective function - is denoted as $J=g(\bm{x};\dot{W})$. It is a function of the state vector $x(t)$, given in (\ref{eqn:probForm}) and wind shear/gradient given by (\ref{eqn:wind}) and (\ref{eqn:logwind}), and is subject to boundary and path constraints \cite{mir2018review}. Only performance indices that are in terms of the energy state of the system are considered here due to their relevance to this work; they are provided in (\ref{eqn:objective}):
\begin{equation}\label{eqn:objective}
    \begin{split}
       { J}&={\max(Energy \: gain)},\\
        {J}&={\max(Total \: energy)}.\\
    \end{split}
\end{equation}

These performance indices (objective functions) in (\ref{eqn:objective}) are consistent with the dynamic soaring literature \cite{gao2017dubins,mir2018review,zhao2004optimal,kim2019deep_total_energy,montella2014reinforcement_energydiff,perez2020neuro_energy,liu2021energy}. The rationale for why the objective functions in (\ref{eqn:objective}) are considered is that the dynamic soaring cycle is characterized as optimal when the system achieves optimal state of energy neutrality (near neutrality); this will be achieved if the system is gaining the most energy it is able to (the first objective function). Equivalently, the system achieves an optimal state of energy if it does not lose the amount of energy it possesses at the beginning of the cycle, i.e., it maximizes the total energy, which cannot exceed the energy level at the beginning of the cycle (the second objective function). It is important to note that the performance indices (objective functions) used in this work are meant to study one cycle of dynamic soaring, i.e., short term behavior of the bird similar to \cite{gao2017dubins,mir2018review,zhao2004optimal,mir2019soaring,mir2021stability,MITbousquet2017optimal,kim2019deep_total_energy,montella2014reinforcement_energydiff,perez2020neuro_energy,liu2021energy}. However, we note that if the bird has a different global objective, such as migration, then the benefits of harvesting energy through dynamic soaring can trade off against the costs of deviating from the desired direction of progress \cite{kempton2022migration}. In such cases, (\ref{eqn:objective}) might not be the right choice for performance indices.  However, in this work, we focus on the local dynamic soaring cycle behavior and not the global behavior of the bird, and we only consider one cycle of dynamic soaring in sufficient windy conditions like those exploited by albatrosses flying in the southern ocean \cite{kempton2022optimization_southernocean}.

Boundary conditions vary with different modes of dynamic soaring. For the basic dynamic soaring mode \cite{mir2018review}, the boundary constraints are 
\begin{equation}\label{eqn:constraints1}
    [z,V,\gamma,\psi]_{t_f}^T=[z+\Delta z, V,\gamma,\psi]_{t_0}^T,
\end{equation}
where $t_f$ and $t_0$ refer to final and initial times respectively,  $T$ is time for a cycle of dynamic soaring, and $\Delta z$ is the net change in altitude/height upon completing the dynamic soaring cycle. This means that the velocity and orientation of the bird/mimicking-object by the end of the dynamic soaring cycle should be the same as the velocity and orientation at the start of the cycle, respectively. Similarly, path constraints for the states and control during the dynamic soaring maneuver are given in (\ref{eqn:constraints2}).
\begin{equation}\label{eqn:constraints2}
    \begin{split}
        V_{\min}<V<V_{\max}, \psi_{\min}<\psi<\psi_{\max},\\
        \gamma_{\min}<\gamma<\gamma_{\max}, x_{\min}<x<x_{\max},\\
        y_{\min}<y<y_{\max},z_{\min}<z<z_{\max} ,\\
        \phi_{\min}<\phi<\phi_{\max}, C_{L_{\min}}<C_L<C_{L_{\max}}.\\
    \end{split}
\end{equation}
These constraints provide the upper and lower bounds to the states and control inputs for the system. These bounds are required by most numerical optimal control methods and solvers to obtain the optimal solution. Since the dynamic soaring maneuver generates high accelerations, an additional path constraint of load factor $n$ (ratio between lift and weight) is added as
\begin{equation}\label{eqn:constraints3}
    n=\frac{L}{mg}\le n_{max}.
\end{equation}
The constraint provided in (\ref{eqn:constraints3}) ensures that the load factor of the mimicking-object does not exceed the maximum limit.

For the convenience of the reader, and due to the broad fields of interest in the topic of this paper, we provide table \ref{tab:technical_terms_sec2}, which includes a collection of technical words/terminologies that have been used in this section with a brief explanation to them. This can help especially readers without much background on nonlinear optimal control literature.

\begin{table}[h]
\centering
\resizebox{0.8\textwidth}{!}{
    \begin{tabular}{|>{\arraybackslash}p{5cm}|>{\arraybackslash}p{7cm}|}
		\hline 
		 \textbf{Technical term} & \textbf{ Explanation} \\
		 \hline
		 Optimal control problem &  A type of problem in control theory that aims at finding the most fitting (optimized) control inputs for a dynamical system under study, such that an objective function need to be satisfied.\\
		 \hline 
		Real-time optimization & It is a kind of optimization methods where the optimal values are calculated instantaneously with time. \\
		 \hline
	\end{tabular}
	}
	\caption{Explanation of some selected control-related technical terminologies used in section 2.}
	\label{tab:technical_terms_sec2}
\end{table}

\section{Extremum seeking systems}\label{ESC}

Recall the motivation and contribution paragraphs from the introduction, we aim at proposing a novel extremum seeking system modeling and control framework for the dynamic soaring phenomenon. Keeping in mind the different communities of interest in the dynamic soaring problem, we prefer to provide in this section some background on extremum seeking systems, and their control structures we mainly use in the following section. Hence, if the reader is familiar with extremum seeking control systems, they may skip this section into the next one. Extremum seeking control system is a model-free, adaptive control technique that stabilizes a dynamical system around the extremum point of an objective function \cite{TanHistory2010,ariyur2003real}. It is particularly useful for systems that lack a proper model or an explicit knowledge of input-output characteristics. It has found application in various disciplines such as formation flight \cite{Chicka2006}, autonomous vehicles \cite{biyik2007autonomous,Cochran2009autonomous}, solar array optimization \cite{Brunton2009,Levya2006}, anti-lock braking system \cite{absDrakunov1995,absDincmen2014}, mobile robots \cite{mobile_applications}, bio-mimicry of fish source seeking \cite{fish1,fish2}, among many others.

\begin{figure}[ht]
    \centering
\includegraphics[width=0.9\textwidth]{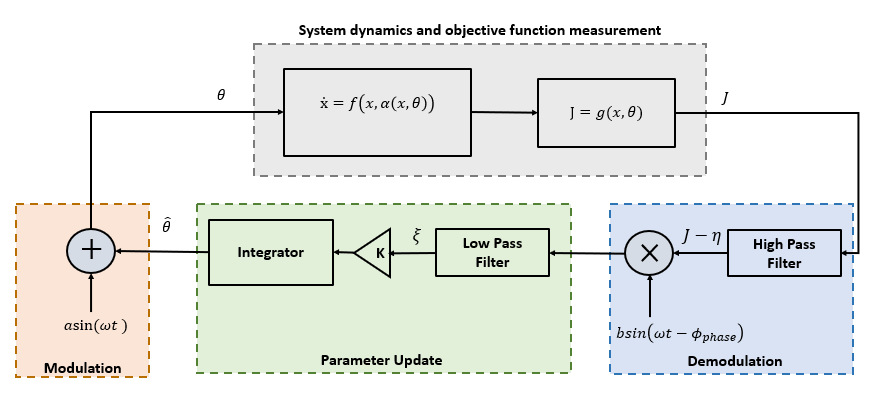}
\caption{Classic structure/design of extremum seeking systems. In this structure, a perturbation signal $a \sin(\omega t)$ is first added to a nominal value of the parameter $\hat{\theta}$ which is the extremum seeking system estimation for the parameter $\theta$. It is initially guessed at the beginning of the problem and then is automatically updated. The system is automatically derived into the optimal value of ($\theta^*$) which in return derives the objective function into its extremum (maximum/minimum). Furthermore, the generated $\theta$ is fed into the system dynamics, which then provides a new evaluation (measurement) of the objective function $J$. Then, the evaluated measurement of $J$ is demodulated using the perturbation signal $b \sin(\omega t-\phi_{phase})$. Finally, the demodulated signal is integrated which then updates the nominal parameter estimation $\hat{\theta}$. The states $\xi$ and $\eta$ are the intermediate states in the structure.}
\label{fig:ESC_Classic}
\end{figure}

We now present two single-input single-output (SISO) structures with the most relevance to this paper. The first structure is the most common extremum seeking structure found in literature, which is usually referred to as the classic structure \cite{KRSTICMain}. This structure can be used for a general nonlinear dynamical system subject to, ideally, a static objective function (its maximum/minimum is fixed and does not change with time, e.g., fixed coordinates). The second structure is an augmented version of the classic structure \cite{krstic2000performance,ESC2_2} and was used, for example, with applications to chemical and industrial process plants. This augmented extremum seeking system can be used for nonlinear dynamic systems subject to a dynamic objective function (its maximum/minimum changes with time, e.g., changing position of a particle). The design of the classic extremum seeking system structure is presented in figure \ref{fig:ESC_Classic}. The parameter $\theta$ is what the extremum seeking controller is updating to reach the extremum of the objective function $J=g({x};\theta)$. Sometimes the parameter $\theta$ is called the maximizer or the minimizer parameter based on the nature of the problem (maximization/minimization). The extremum seeking classic structure consists of three main steps. In the first step, a perturbation signal (modulation) -- usually a sinusoidal signal in the form $a \sin(\omega t)$, where $a$ and $\omega$ are the amplitude and frequency of the signal -- is added to a nominal value of the parameter $\hat{\theta}$, which is an estimation for $\theta$. Initially, at the beginning of the problem, we start with an initial guess value of it. The generated $\theta$ is fed into the system dynamics $\dot x = f(x,\theta)$, which then provides a new evaluation (measurement) of the objective function $J$; we note here that the mathematical expression of the objective function is not required, only the measurement. 
% Make a note that a more complex form can be considered if we have a particular representation of the parameter theta inside f , alpha. But we do not consider it here for simplicity and focus of the paper.
In the second step, the evaluated measurement of $J$ is demodulated using a perturbation signal $b \sin(\omega t-\phi_{phase})$ of the same frequency but with amplitude $b$ and phase lag of $\phi_{phase}$. The result of that determines if the estimation of $\theta$ needs to increase or decrease. In the third and final step, the demodulated value is integrated which then updates the nominal parameter estimation $\hat{\theta}$. The high pass filter removes the ``dc component" in the signal before the demodulation step and the low pass filter extracts the ``dc component" for parameter update. Note that the ``dc component" is the mean amplitude of the signal. The use of these filters is optional. The reader can refer to \cite{KRSTICMain} for more technical explanation and the Mathworks web page \cite{matlab_esc} for a much more graphical and easier explanation. 

\begin{figure}
    \centering
    \includegraphics[width=0.8\textwidth]{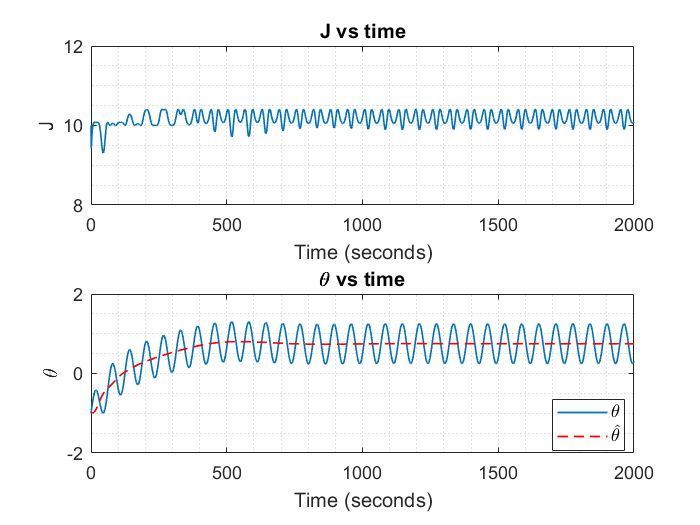}
    \caption{An example to demonstrate the utility of the classic extremum seeking control. The upper plot shows the evolution of the objective function over time, and it oscillates around (stabilizes about) the extremum (maximum in this case) value of $J^* =10.41$. Similarly, the bottom plot provides the evolution of parameters $\theta$ and $\hat{\theta}$ overtime. Result shows that the estimate $\hat{\theta}$ approaches the optimal value of $\theta^*=0.88$ and the parameter $\theta$ oscillates/stabilizes about it.}
    \label{fig:esc_example}
\end{figure}

To demonstrate the utility of the classic extremum seeking control, we will adopt a maximization problem from \cite{tan2009global}. Let us consider a dynamical system with the following two dynamical equations
\begin{equation}\label{eqn:esc_example_eqn}
\dot{x}_1 = -x_1+x_2;\:
\dot{x}_2 = x_2+u.\:
\end{equation}
Let the control input be $u=-x_1-4x_2+\theta$. Ideally, as explained earlier, we do not need a mathematical expression for the objective function, however, for our aimed numerical illustration, we assume access to the expression/formula of the objective function, $J(\theta)= -\theta^4+\frac{8}{15} \theta^3+ \frac{5}{6} \theta^2 +10$. This static objective function has a global maximum of $J^* = 10.41$ corresponding to $\theta^* = 0.88$. Now, we implement the above problem using the extremum seeking structure in figure \ref{fig:ESC_Classic}, with the system dynamics provided in (\ref{eqn:esc_example_eqn}). The amplitude and frequency of modulating and demodulating signals are taken as $a = 0.5$ and $\omega = 0.1 $. The high pass filter and low pass filter have the form of $s/(s+\omega_h)$ and $w_l/(s+\omega_l)$ respectively, with $\omega_h =0.03$ and $\omega_l = 0.01$. The initial states of the system are taken as $[x_1,x_2]=[0,0]$ with the initial estimation value of the parameter as $\hat{\theta} = -1$. We ran the simulation for 2000 seconds and the result is shown in figure \ref{fig:esc_example}. The result shows the extremum seeking system's ability to derive the parameter $\theta$ in a neighborhood (stabilizes about) the optimal parameter value $\theta^*$. This has caused the output (the objective function) $J$ to stabilize about its optimal (maximum) value $J^*$ which corresponds to $\theta^*$.

\begin{figure}[ht]
    \centering
\includegraphics[width=0.75\textwidth]{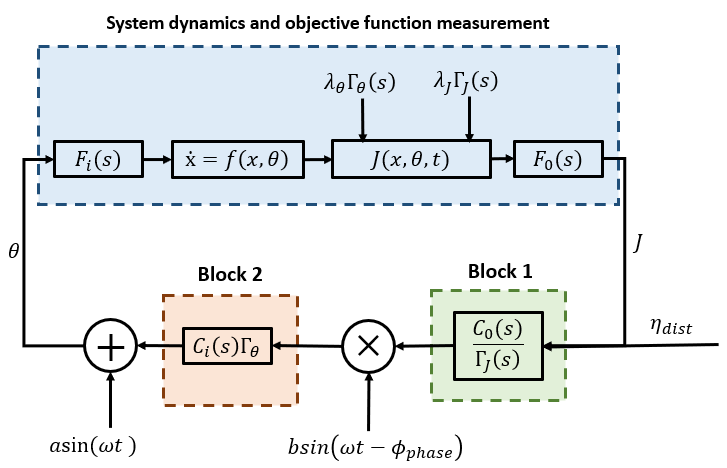}
\caption{Design of an augmented extremum seeking structure used for extended nonlinear dynamical systems. Here, $F_i(s)$ and $F_0(s)$ denote an additional (optional) input and an additional (optional) output functions that can extend the dynamics based on the application. The parameter $\theta$ is the maximizer/minimizer parameter of interest as is the case with the classic structure in figure \ref{fig:ESC_Classic}. The objective function $J$ is the output and it can be depending on time as well as $\theta$. $\lambda_\theta \Gamma_\theta(s)$ and $\lambda_J \Gamma_J(s)$ denote the Laplace transform of the time-varying optimal input $\theta^*(t)$ and output $J^*(t)$ as shown in (\ref{eqn:laplace_transform}). $C_0 (s)$ and $C_i(s)$ are the compensators, design of which depend on the dynamics $F_0(s)$ and $F_i(s)$, as well as the reference signals $\Gamma_\theta(s)$ and $\Gamma_J(s)$. The design procedure is provided in appendix \ref{sec:app_algorithm}. Finally, $\eta_{dist}$ represents a general disturbance that is reflected in the measurement of objective function due to disturbances/perturbations in the computations or noises in the sensors themselves.}
\label{fig:ESCPlant}
\end{figure}

Now, if the extremum of the objective function is changing with time, then an augmented structure of the extremum seeking classic structure can be used \cite{krstic2000performance,ESC2_2} as mentioned earlier; this augmented classic extremum seeking structure/design is shown in figure \ref{fig:ESCPlant}.
 Here, $F_i(s)$ and $F_0(s)$ denote optional external input and output dynamics that can be used for better representation of the nonlinear system dynamics in some applications. Let $\theta^*(t)$ and $J^*(t)$ be the time-varying optimal input and optimal output/objective-function, respectively. Then their Laplace transform can be denoted as
\begin{equation}\label{eqn:laplace_transform}
    \begin{split}
        L[\theta^*(t)]&=\lambda_\theta \Gamma_\theta(s)\\
        L[J^*(t)]&= \lambda_J \Gamma_J(s).
    \end{split}
\end{equation}
Having $\theta^*(t)$ and $J^*(t)$ as time-varying, allows for the possibility of having to optimize a system whose commanded operation is not a constant, hence, time-varying maxima or minima can be tracked. $C_0 (s)$ and $C_i(s)$ are the compensators, design of which depend on the dynamics $F_0(s)$ and $F_i(s)$, as well as reference signals $\Gamma_\theta(s)$ and $\Gamma_J(s)$ (see appendix \ref{sec:app_algorithm}). For example, if the input signal we are tracking has the nature of a ramp signal ($\theta$ is proportional to time), $\theta^*(t)\propto t$, then 
% $\Gamma_\theta(s) = 1/{s^2}$,
\begin{equation}\label{eqn:laplace}
    L[\theta^*(t)]=\lambda_\theta \Gamma_\theta(s) = \frac{\lambda_\theta}{s^2},
\end{equation}
which means we will get a double integrator in the feedback loop which increases the possibility of having an unstable system. But, by choosing proper $C_i(s)$, the relative degree of the loop can be reduced. Thus, the compensators $C_0 (s)$ and $C_i(s)$ are important design tools for satisfying stability conditions. 

\begin{figure}
    \centering
    \includegraphics[width=0.8\textwidth]{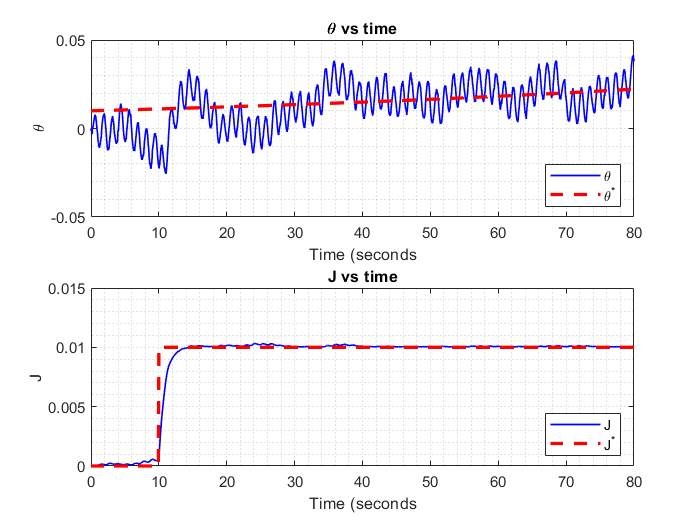}
    \caption{An example to demonstrate the utility of the augmented extremum seeking controller with nonlinear dynamics with dynamic (time changing) objective function. The plot in the top shows the evolution of the input parameter $\theta$ (blue) tracking the time-varying optimal parameter ${\theta}^*$ (red) over time. The bottom plot shows the evolution of time-varying objective function $J$ (blue) tracking its extremum (optimal value) which is also time-varying $J^*$ (red) over time. The optimal parameter is changing with time and the optimal objective function has a step jump and yet the structure is able to track them. }
    \label{fig:esc_example_plant}
\end{figure}

To demonstrate the utility of this structure, we present an example similar to \cite{ariyur2003real}. We take $F_i(s) = {(s-1)}/{(s^2+3s+2)}$, $F_0 = {1}/{(s+1)}$, $J(\theta) = J^*(t)+(\theta-\theta^*(t))^2$, with $J^*(t) = 0.01u(t-10)$, $\theta^*(t) = 0.01e^{0.01t}$, where $u(t-\tau)$ represents the unit step function at time $\tau$. By Laplace transform, we get $\lambda_J \Gamma_J (s)=(0.01e^{-10s})/s$ and $\lambda_{\theta} \Gamma_{\theta}(s) = 0.01/(s-0.01)$. Using the algorithm presented in section \ref{sec:app_algorithm}, the block 1 and block 2 are designed as $s/(s+5)$ and $50(s-4)/(s-0.01)$ respectively, with $\omega =5$ rad/s, $a=0.05$ and $\phi_{phase} = 0.8$. The simulation is run for 80 seconds and the result is shown in figure \ref{fig:esc_example_plant}. The result shows that the optimal parameter and optimal objective function are tracked even when they are changing with time.

Similar to the previous section, we provide table \ref{tab:technical_terms_section3}, which includes a collection of technical words/terminologies that have been used in this section with a brief explanation to them. This can help the readers who do not have much background on extremum seeking control literature.

\begin{table}[ht]
\centering
\resizebox{0.8\textwidth}{!}{
    \begin{tabular}{|>{\arraybackslash}p{3cm}|>{\arraybackslash}p{10cm}|}
		\hline 
		 \textbf{Technical term} & \textbf{ Explanation} \\
		 \hline
		 SISO & It stands for single-input single-output. SISO control systems are the systems which take only one input and provide only one output.\\
		 \hline
		 Modulation & It is a process of modifying a given signal to carry or be assigned some additional information.\\
		 \hline
		  Demodulation& It is a process of extracting the original information from the modulated signal.\\
		  \hline
		 High pass filter & It is a device/method that passes the signals with frequencies higher than a certain cut-off frequency and filters out the rest.\\
		 \hline
		  Low pass filter & It is a device/method that passes the signals with frequencies lower than a certain cut-off frequency and filters out the rest.\\
		 \hline
		 $s$& It represents the Laplace transform variable.\\
		 \hline
		 Optimal input & An input that provides the maximum (or minimum) value of the desired objective function.\\
		 \hline
		 Optimal output & The maximum (or minimum) value of the desired objective function.\\
		 \hline
	\end{tabular}
	}
	\caption{Explanation of some selected control and technical words/terminologies used in section 3.}
	\label{tab:technical_terms_section3}
\end{table}
\section{Extremum seeking characterization of the dynamic soaring phenomenon}\label{ESC_DS}
 In this section, we hypothesize, and provide our proof of concept, that the dynamic soaring phenomenon can be characterized (described, modeled, simulated, mimicked, and controlled) by extremum seeking systems. First, we explain our philosophy behind the linkage between extremum seeking systems and dynamic soaring. Second, we provide our proposed extremum seeking models and control structures for dynamic soaring. Recall from section \ref{ESC} that extremum seeking systems are both models and controls at the same time; they perform both the functions of trajectory planning and tracking as depicted in figure \ref{fig:traj_plan}. Third, we provide simulation results to illustrate the effectiveness of the proposed framework and to support further our hypothesis; our simulations also include comparisons with powerful optimal control methods found in the literature. Lastly, we provide a comprehensive discussion on the proposed approach and how it compares with literature, its stability, applicability, potentials as well as limitations.

\subsection{Philosophy behind characterizing dynamic soaring as an extremum seeking system}
\begin{figure}[ht]
    \centering
\includegraphics[width=0.8\textwidth]{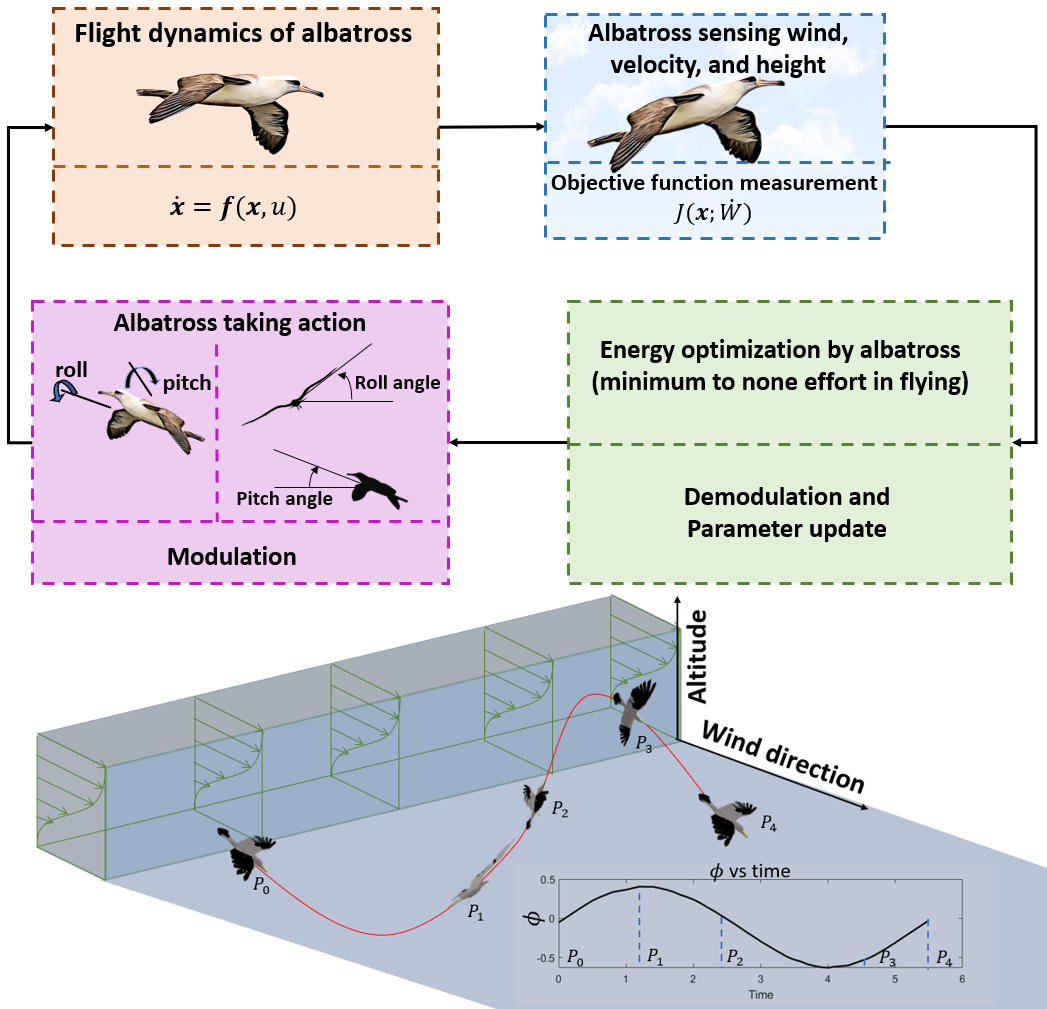}
\caption{Parallelism between dynamic soaring and extremum seeking systems. The (pitching and rolling) variation/perturbation action taken by the bird can be considered as the modulation step. These actions then affect the flight dynamics state space including the velocity and orientation (see (\ref{eqn:dynamics})). Then, sensing the change in the wind, velocity and orientation is similar to measuring the objective function (see (\ref{eqn:objective})) that extremum seeking systems do not require its mathematical expression. Finally, the assessment of the measurement determines if the energy conservation is met properly or not. Hence, updating the action as per the feedback is similar to the demodulation and parameter update steps. The top four blocks depict the process and the parallelism hypothesized above. The bottom two wide graphs explain further how the action (control) of the bird is done autonomously during the dynamic soaring cycle. For instance, let us track only the rolling action of the bird while it is tracking the red curve (the three-dimensional) dynamic soaring trajectory. At the point $P_0$, the bird starts the maneuver with 0 degrees rolling angle. Based on the action$\rightarrow$sensing$\rightarrow$feedback process described above, it changes/increases slowly its rolling consistently until reaches about the point $P_1$. If it continues changing rolling as it was (by increasing rolling), it would sense less energy conservation, thus, it decreases its rolling until it gets back to 0 degrees again, point $P_2$. This time, continuous change in the negative rolling direction is generating a plausible assessment of energy conservation until it reaches about the point $P_3$, where the bird then has to change course, increasing its rolling leading to the point $P_4$. The rolling action is provided in radian vs. time in the right bottom graph.}
\label{fig:philosophy}
\end{figure}

The dynamic soaring phenomenon clearly asserts that albatrosses travel hundreds of miles while spending a low amount of energy \cite{mir2018review}. 
% This implies that they should have a mechanism and process that enables them to find a path that maximizes the energy gained from the wind so that they minimize the effort (energy spending) in flying. 
This implies that they should have a mechanism to harvest energy from other means; this is shown -- with experimental validation -- to be done via finding a flight path that enables harvesting energy from the wind \cite{sachs2013experimental}.
It is also important to note that they find said path autonomously and in real-time. A natural way to perform such a maneuver -- assuming the presence of sufficient wind condition  -- is to take periodic variation/perturbation in possible actions (controls), e.g., pitching and/or rolling, and then by sensing/measuring the corresponding changes in wind, velocity and height, the bird determines its energy state and whether it needs to change actions to spend less energy; recall that the nostrils of the albatross act as an airspeed sensor \cite{pennycuick2008information_nostrils,brooke2002gusts_nostrils}. Then, the actions are updated accordingly and the process is autonomously performed until the path/trajectory is both identified and tracked. The periodic nature of the perturbations can be easily memorized and repeated. The above-mentioned hypothesis matches that of extremum seeking systems as explained in more detail earlier in section \ref{ESC}. Figure \ref{fig:philosophy} draws the hypothesized parallelism between the dynamic soaring phenomenon and extremum seeking systems. By recalling the blocks and steps representing the functionality of extremum seeking systems in figure \ref{fig:ESC_Classic} from the previous section, and comparing it with the parallelism provided in figure \ref{fig:philosophy}, the following linkage becomes clearer.
The pitching and rolling variation/perturbation action taken by the bird can be considered as the modulation step in extremum seeking. Then, sensing the change in the wind, velocity and height is similar to measuring the objective function that extremum seeking systems do not require its mathematical expression. Finally, updating the action as per the energy assessment feedback is similar to the demodulation and parameter update step. Consequently, we hypothesize that extremum seeking systems can be considered a natural and comprehensive modeling and control system framework for dynamic soaring. In the next subsection, we provide our proposed extremum seeking model for dynamic soaring along with the proposed control structures for the problem.

\subsection{Proposed extremum seeking modeling and control structures for dynamic soaring}\label{ESC_DS_structure}
As mentioned in section \ref{problem_setup}, the actions taken by the bird/mimicking-object for the dynamic soaring system are typically the control inputs $\bm{u}(t)=[C_L,\phi]$ in (\ref{eqn:probForm}), representing the pitching and rolling actions, respectively. However, the maneuver can be performed with a free (varying) control of $C_L$ as well as a constant $C_L$ as shown in \cite{sachs2005minimum}. We choose the case with the constant $C_L$ to make the dynamic soaring problem a single-input single-output (SISO) system with the rolling $\phi$ as the free (varying) single control input, and the objective function $J$ in (\ref{eqn:objective}) being the single output of the system. This allows us to take advantage of the mature and developed tools for single-input-single-output extremum seeking controllers, available in the extremum seeking literature \cite{ariyur2003real,scheinker2017model}. It is also our assessment that since extremum seeking application to the dynamic soaring problem is novel and without precedents, the introductory results of this problem are better to be in a single-input single-output basis before increasing the complexity further.

%  Since dynamic soaring can be seen as a system with a relatively slow frequency that represents the flight dynamics of a fairly large bird (albatross) and it takes a relatively long time to perform its cycle, we are also going to implement the problem using the classic structure in figure \ref{fig:ESC_Classic} which can work for some time-varying outputs, especially when the system dynamics is slow. This is especially true if the objective function is the maximum of the total energy as the maximum here is constant: the amount of energy the system possesses at the beginning of the cycle. This represents then a static (time-invariant) objective function, which fits the characteristics of the classic structure \cite{KRSTICMain} in figure \ref{fig:ESC_Classic}.

Recall the classic extremum seeking system explained earlier in section \ref{ESC}. Now, we will customize that structure to propose it as a novel modeling and control framework for dynamic soaring. First, we may maintain the albatross/mimicking-object flight dynamics model given in (\ref{eqn:dynamics}). The single input (control) is the rolling action $\phi$ -- recall also figure \ref{fig:philosophy}. The output, which is the objective function $J$ is assumed to depend generally on the state variables of the flight dynamic system and the wind shear/gradient, hence $J=g(\bm{x},\dot{W})$. We set the parameter $\theta$, the maximizer or minimizer, of the objective function $J$, as $\theta=\phi$. That is,  the action/control-input of the rolling will be the parameter maximizing or minimizing the objective function of the system. This sets up the extremum seeking system to steer the albatross/mimicking-object system, utilizing the optimal rolling action/control, towards the extremum (maximum/minimum) of the objective function. The high pass filter and low pass filter have the form of $s/(s+\omega_h)$ and $w_l/(s+\omega_l)$ respectively. The modulating and demodulating signals are $a \sin(\omega t)$ and $b \sin(\omega t-\phi_{phase})$. Recall from section \ref{ESC} that $\omega$, $\omega_h$, and $\omega_l$ are the frequencies corresponding to the input/modulating signal, high pass filter, and low pass filter, respectively. Similarly, $\phi_{phase}$ is the phase lag between the modulating and demodulating signals.

\begin{figure}[ht]
    \centering
\includegraphics[width=1\textwidth]{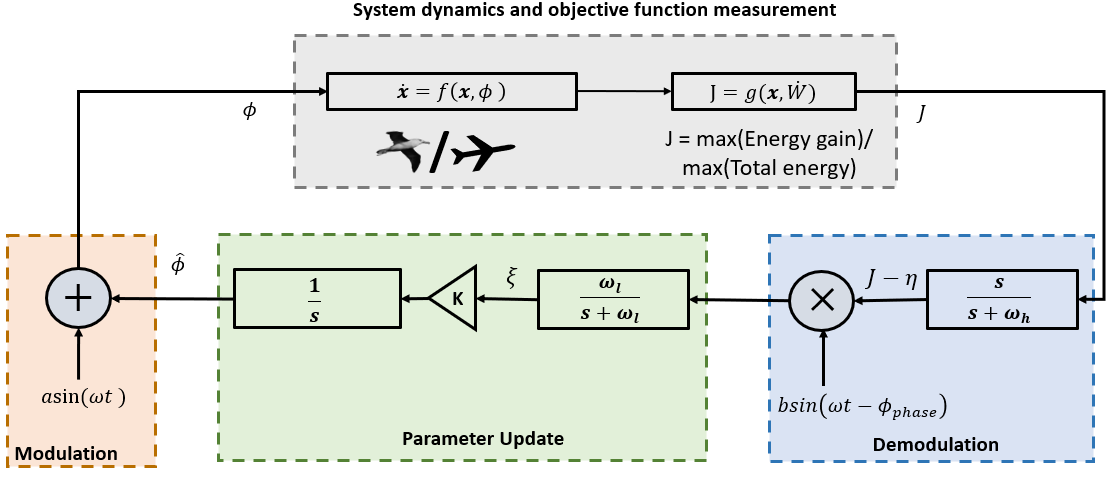}
\caption{Customized classic extremum seeking structure for the dynamic soaring. The control input to the flight dynamic system in (\ref{eqn:dynamics}) is the rolling angle ($\phi$). The symbols and blocks are consistent with the symbols and blocks described in figure \ref{fig:ESC_Classic}.}
\label{fig:ESC_Classic_custom}
\end{figure}

We propose the classic extremum seeking system in figure \ref{fig:ESC_Classic_custom} as a modeling and control system for dynamic soaring.
Equivalently, the proposed time domain mathematical model and control system of the dynamic soaring problem is provided in \eqref{eqn:esc_x_dynamics}- \eqref{eqn:esc_eta_dynamics}:
\begin{align}
    \dot{\bm{x}} &= f(\bm{x},u), \label{eqn:esc_x_dynamics}\\
    \dot{\hat{\phi}} &= K \xi,\\
    \dot{\xi} &= -\omega_l \xi +\omega_l(J-\eta) b \sin (\omega t -\phi_{phase}), \label{eqn:esc_xi_dynamics}\\
    \dot{\eta} &= -\omega_h \eta+ \omega_h J, \label{eqn:esc_eta_dynamics}
\end{align}
 where $\dot{\bm{x}}=f(\bm{x},u)$ in (\ref{eqn:esc_x_dynamics}) is the flight dynamics model of an albatross/mimicking-object as in (\ref{eqn:dynamics}) with a rolling control input $u=\phi=\hat{\phi}+a \sin (\omega t)$, where $\hat{\phi}$ is the estimated value of the control input $\phi$. Similarly, $a$ and $\omega$ are the amplitude and the frequency of the input signal, $K$ is some positive constant as long as the problem is a maximization one (negative if the problem is a minimization one), and $\xi$ and $\eta$ are intermediate states in the structure as shown in figure \ref{fig:ESC_Classic_custom}.

 As discussed in section \ref{ESC}, the objective function $J$ in the model above
 (see equations (\ref{eqn:esc_xi_dynamics})-(\ref{eqn:esc_eta_dynamics})) is better to be assumed static (not time dependant) for better applicability of the classic extremum seeking structure. This will not be an issue with the objective function aiming at maximizing the total energy (\ref{eqn:objective}), as the total energy remains constant or near-constant during the dynamic soaring cycle with a fixed maximum, i.e., maximizing the total energy is a static objective function. However, the energy gain has been observed to vary with time \cite{sachs2005minimum}, i.e., an objective function aiming at maximizing energy gain may be dynamic (its maximum changes with time) and not static. As we will show later in section \ref{simulationSection} in our simulations, the classic structure of extremum seeking can still incorporate time-varying objective functions as long as the system is slow (the case here as the dynamic soaring cycle does not possess high frequencies with a relatively large periodic cycle). However, we believe it might be necessary to introduce also a structure that is well-defined and able to tolerate time-changing objective functions such as the augmented extremum seeking system explained in section \ref{ESC}. Our reasoning for choosing said augmented extremum seeking structure can be summarized in three technical points: (i) it allows for fast-tracking of the output even in cases with fast changes (recall the solved example in section \ref{ESC}, figure \ref{fig:esc_example_plant}), which is arguably what the bird is doing to avoid deflecting away from its optimized flight path, (ii) it guarantees stability if certain assumptions and conditions (verifiable and can be constructed by design) are satisfied unlike many other extremum seeking structures where the stability property is characterized by general theorems and conditions that are not customized to the particular problem; and (iii) it allows the optimal output (objective function) to be time-varying, which matches the nature of the dynamic soaring problem that may maximizes or minimizes its output as a function of time. Next, we customize the augmented extremum seeking structure to adopt the dynamic soaring problem.

\begin{figure}[ht]
    \centering
\includegraphics[width=1\textwidth]{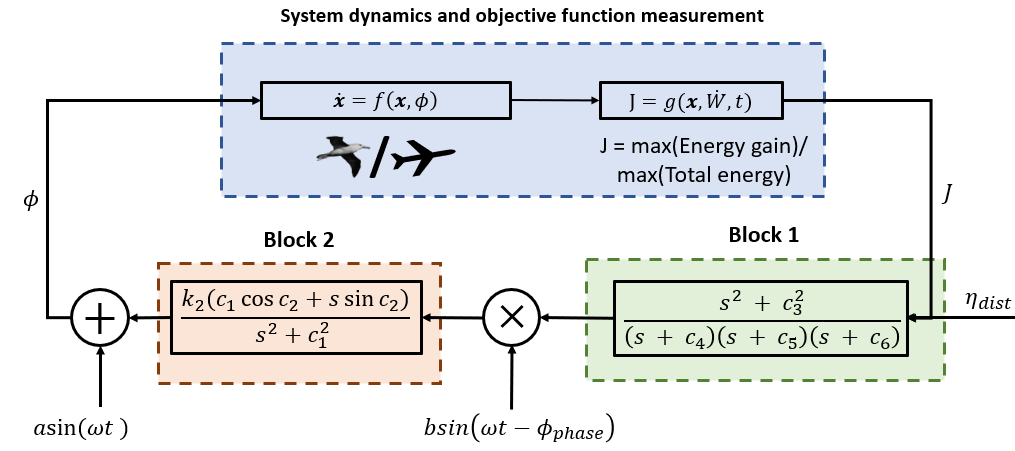}
\caption{Customized augmented extremum seeking structure for the dynamic soaring problem, especially applicable if the objective function is time-varying with a time-varying maximum or minimum. The control input to the flight dynamics system in (\ref{eqn:dynamics}) is the rolling angle $\phi$. The symbols and blocks are consistent with the symbols and blocks in figure \ref{fig:ESCPlant} with no external input/output dynamics: $F_i(s)=1$ and $F_0(s)=1$.}
\label{fig:ESCPlant_custom}
\end{figure}

Similar to the classic structure, we maintain the albatross/mimicking-object flight dynamics model given in (\ref{eqn:dynamics}). The single input (control) is the rolling $\phi$. The output is $J=g(\bm{x},\dot{W})$. Similarly, we set the maximizer or minimizer parameter $\theta$ as $\theta=\phi$. Recall the augmented extremum seeking structure in figure \ref{fig:ESCPlant}. We note that for the dynamic soaring problem, there are no external dynamics cascaded with the albatross/mimicking-object model, i.e., there is no $F_i(s)$ or $F_0(s)$. In other words, there is no input dynamics $F_i(s)$ between the control input and the flight dynamics model as the control input affects directly the flight dynamic system without any additional transfer functions. Similarly, there is no output dynamics $F_0(s)$ between the flight dynamics model and the measurement of the objective function as $J$ is measured based on the flight dynamic state variables and the wind shear without any transfer function affecting the measurement. Consequently, for the augmented extremum seeking structure, we set $F_i(s)=1$ and $F_0(s)=1$. With the proper design choices of Block 1 and Block 2, we propose the augmented extremum seeking structure in figure \ref{fig:ESCPlant_custom} as a novel control system for dynamic soaring.
% We denote the optimal input and output as $\phi^*$ and $J^*$ respectively. Now, their Laplace transforms are denoted as
% \begin{equation}
%     \begin{split}
%         L[\phi^*(t)]&=\lambda_\phi \Gamma_\phi(s)\\
%         L[J^*(t)]&= \lambda_J \Gamma_J(s),
%     \end{split}
% \end{equation}
% where $\lambda_\phi$ and $\lambda_J$ are constants. If $\phi^*$ and $\theta^*$ are constants (step functions), then the structure in figure \ref{fig:ESCPlant} is equivalent, and reduces, to the classic structure \cite{KRSTICMain} in figure \ref{fig:ESC_Classic}.

Now, we return back to the design of the augmented extremum seeking structure for dynamic soaring. We here provide a brief explanation for our steps based on the guidelines and design algorithm  provided as algorithm 2.1 in \cite{ESC2_2} -- also provided in appendix \ref{sec:app_algorithm} of this paper. However, for more technical details, we recommend the reader to read appendix \ref{sec:app_algorithm}. First, we need to find $\Gamma_J$ and $\Gamma_\phi$ corresponding to Block 1 and Block 2 respectively. For that, we make a guess on the nature of $\phi^*(t)$ and the nature of $J^*(t)$, the optimal input and output of the system.
% Following the design steps of Algorithm 2.1 in \cite{ESC2_2}, we need to make an assumption on the nature of $\theta^*(t)$ and the nature of $J^*(t)$. These assumptions are not the ``equivalent" expression for $\theta^*(t)$ and $J^*(t)$. As a matter of fact, neither $\theta^*(t)$ nor even $J^*(t)$ need to be known a priori.
From many observations made on dynamic soaring maneuvers, theoretically as in \cite{mir2018optimal,gao2017dubins,mir2019soaring,MITbousquet2017optimal} and experimentally as in \cite{sachs2013experimental,yonehara2016flight}, one can deduce that the nature/guess of $\phi^*(t)$ and $J^*(t)$ is a periodic or sinusoidal-like behavior. A further explanation of this reasoning is also provided in appendix \ref{sec:app_algorithm}. Hence, we assume the optimal input and output nature/guess as sinusoidal signals with phase shift, i.e., $\phi^*(t)=\sin(c_1 t+c_2)$ and $J^*(t)=\sin(c_3t+c_0)$, where $c_1,c_2,c_3,c_0$ are some parameters. By using Laplace transform we find $\Gamma_\phi$ and $\Gamma_J$ as
\begin{equation}\label{eqn:gammas}
    \begin{split}
        \Gamma_\phi&=\frac{c_1\cos c_2 + s \sin c_2}{c_1^2+s^2},\\
        \Gamma_J&= \frac{c_3 \cos c_0 + s \sin c_0}{c_3^2+s^2}.
    \end{split}
\end{equation}
Next, we need to design the compensators $C_0$ and $C_i$ in Block 1 and Block 2 in figure \ref{fig:ESCPlant_custom}. Since we have taken $F_i=1$ and $F_0=1$, the compensators $C_i$ and $C_0$ only depend on $\Gamma_\phi$, $\Gamma_J$, the nature/guess of optimal parameter $\phi^*$ and optimal objective value $J^*$. This means that the compensators are designed in such a way that this structure of extremum seeking is able to utilize the information about the nature/guess of the optimal input and output to better track the actual optimal output. Hence, this structure augments the classic extremum seeking structure.
By following the algorithm in appendix \ref{sec:app_algorithm}, we get $C_i = 1$ and $C_0 =(c_3 \cos c_0 + s \sin c_0)/((s+c_4)(s+c_5)(s+c_6))$, where $c_4,c_5$ and $c_6$ are some constants. Thus, we get $C_0/\Gamma_J = (s^2+c_3^2)/((s+c_4)(s+c_5)(s+c_6))$ and $C_i\Gamma_{\phi} = (c_1\cos c_2 + s \sin c_2)/(c_1^2+s^2)$.
This completes the customization and further explains the origin of the transfer functions and blocks in figure \ref{fig:ESCPlant_custom}.

\subsection{Simulations and comparative results}\label{simulationSection}
In this subsection, we perform simulations for the proposed extremum seeking systems customized for dynamic soaring in the previous subsection. Moreover, we compare said simulation results with the solutions obtained from an optimal control numerical optimizer. The first structure of extremum seeking corresponding to figure \ref{fig:ESC_Classic_custom}-- also provided in (\ref{eqn:esc_x_dynamics})-(\ref{eqn:esc_eta_dynamics}) is referred to as ``ESC1" throughout the simulations. The second structure corresponding to figure \ref{fig:ESCPlant_custom} is referred to as ``ESC2".
 As for the numerical optimizer, we have used GPOPS2 \cite{gpops2}, which is an optimization software compatible with MATLAB\textsuperscript{\textregistered} and uses hp-adaptive Gaussian quadrature collocation methods and sparse nonlinear programming.
To broaden our results and proof of concept, we will perform simulations involving two different performance indices (objective functions) and two different wind shear models.  The chosen wind shear models are the logistic and logarithmic models, as described in section \ref{problem_setup} with parameters provided in table \ref{tab:Windparams}, whereas the chosen performance indices (objective functions) are provided in (\ref{eqn:objective}).

Next, we derive the expressions of the chosen performance indices. The expression corresponding to specific total energy as a performance index is
\begin{equation}\label{eqn:specificTE}
    e = TE/mg =z+\frac{V^2}{2g},
\end{equation}
where $TE$ represents the total energy, $e$ represents specific total energy, $z$ is the altitude, $V$ is the velocity, $m$ is the mass of the bird/mimicking-object, $g$ is the acceleration due to gravity. Now, for the expression corresponding to specific energy gain, taking derivative of (\ref{eqn:specificTE}) with respect to time and using the equations of $\dot{V}$ and $\dot{z}$ in (\ref{eqn:dynamics}), we get
\begin{align}
    \dot{e}&=\dot{z}+\frac{V\dot{V}}{g}\\
    &= Vsin\gamma + \frac{V}{g}[-D-g\sin\gamma -\dot{W}\cos \gamma \sin \psi] \\
    &=-\frac{DV}{g}-\frac{V\dot{W}\cos \gamma \sin \psi}{g}.\label{eqn:energyGain}
\end{align}
According to (\ref{eqn:energyGain}), the term $-V\dot{W} \cos \gamma \sin \psi /g $ is the one determining energy gain from the wind and is taken as a performance index. Thus, we use $J_1 = max(-(V\dot{W}\cos \gamma \sin \psi)/g$)  as our first performance index (objective function maximizing energy gain) and $ J_2=max(z+\frac{V^2}{2g})$ as the second performance index (objective function maximizing total energy).
% We use this performance index measure for both extremum seeking structures and GPOPS2, whereas for the DC method of optimization, we do not use the performance index measure and solve it as the ``constraints satisfaction" problem as done in \cite{MITbousquet2017optimal}. For some optimal control problems in dynamic soaring literature, the problem is so heavily constrained that it does not require an objective function. This is a problem, we hope this work will help to resolve. \textit{It is important to emphasize that all bounds and constraints in \eqref{eqn:constraints1}, \eqref{eqn:constraints2} and \eqref{eqn:constraints3} are not applied or required in ESC1 and ESC2, however, they are a must for GPOPS2 and DC method.}
Next, we explain each extremum seeking structure used in the simulations. For ESC1, the high pass filter has the expression $s/(s+h)$ with $K=k_1$ and the integrator as $1/s$. We have omitted the use of the low pass filter as it is optional. The modulating and demodulating signals are in the form $a_1 \sin(\omega_1 t)$ and $b_1 \sin(\omega_1 t + \phi_{phase_1})$ respectively. For ESC2, we design the Block 1 in the form of $ (s^2+c_3^2)/((s+c_4)(s+c_5)(s+c_6))$ and Block 2 as $k_2(s \sin(c_2)+c_1\cos(c_2))/(s^2+c_1^2)$ as exactly provided in figure \ref{fig:ESCPlant_custom}. The modulating signal is in the form $a_2 \sin(\omega_2 t)$, whereas the demodulating signal is in the form $b_2 \sin(\omega_2 t + \phi_{phase_2})$.  For all simulations, we apply $C_L=1.5$. The parameters of the albatross flight dynamics are taken from table \ref{tab:params} .
% The initial state of the simulation is taken as $[x,y,z,V,\gamma,\psi]=[-16,15,0,14,-0.4,-1.16]$ and it is fixed for all simulations.

\begin{figure*}[ht]
\centering
	\includegraphics[width=1\textwidth]{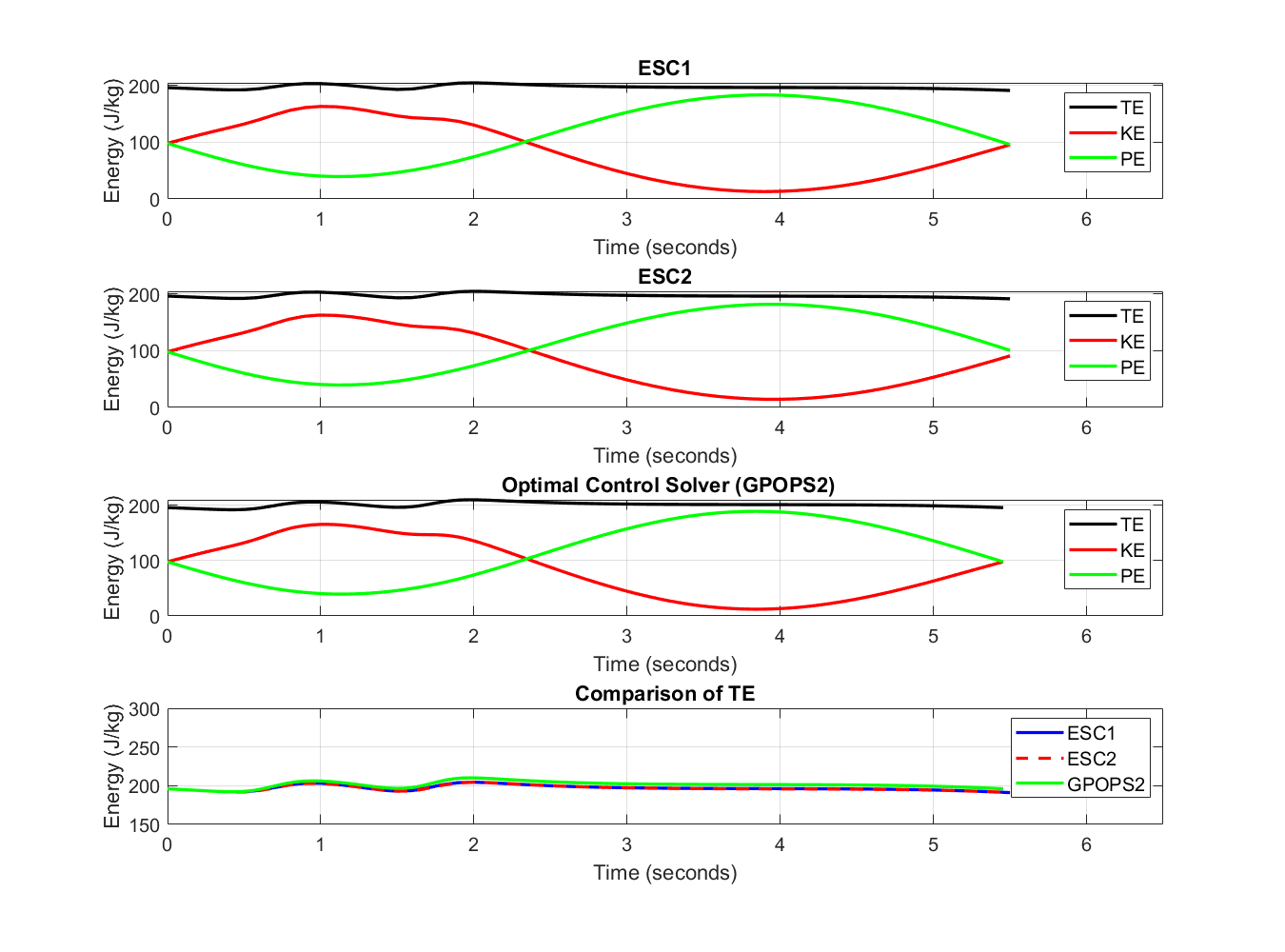}
	\caption{Comparison of total energy during a cycle of dynamic soaring obtained using ESC1, ESC2, GPOPS2 with logistic wind profile and $J_1$ as the performance index (case 1). The results obtained from all three methods are similar. Also, the total energy is near-constant in all of the three methods which indicates the success of our implementation.}
	\label{fig:case1_allenergies}
\end{figure*}

\begin{figure*}[ht]
\centering
	\includegraphics[width=1\textwidth]{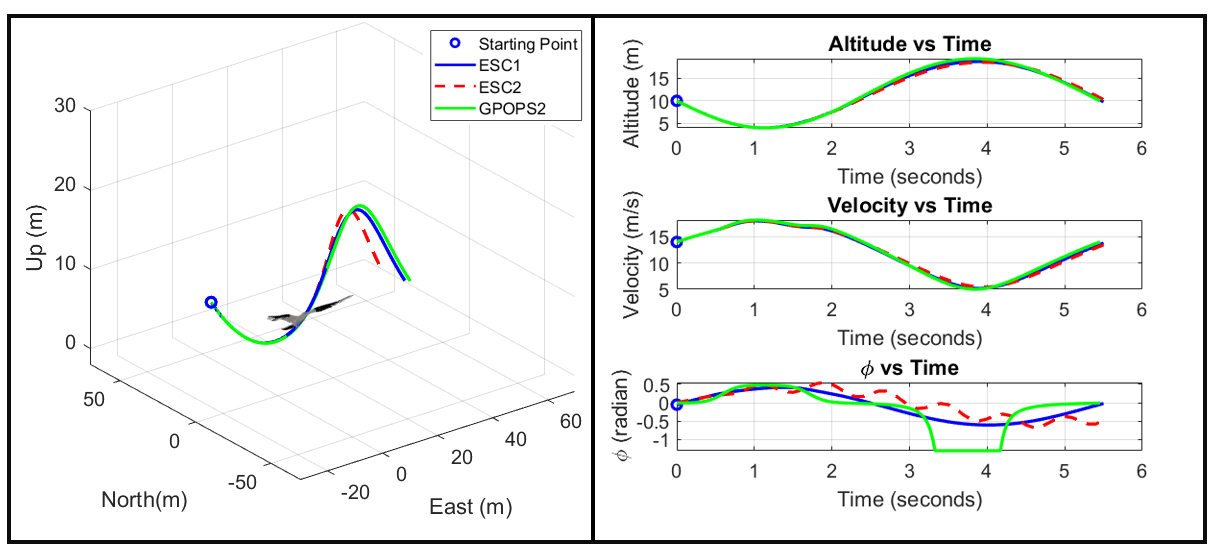}
	\caption{ Case 1 of our simulation results. It provide a comparison of three dimensional trajectories, velocity, altitude and control input $\phi$ of an albatross/mimicking-object conducting dynamic soaring. In the simulations, we used two extremum seeking systems ESC1 and ESC2, and a powerful nonlinear optimal control solver, GPOPS2. The wind conditions in the simulations  were considered to follow logistic wind profile. The objective function $J_1$, the performance index, here is aiming at maximizing the energy gain. The trajectories obtained with different methods are similar to each other. Similarly, the evolution of the parameters with time is comparable even though they are obtained from three different methods. The 3D trajectories are provided in the left. The altitude, velocity, and control input (rolling) are provided in the right part of the figure.}
	\label{fig:case1_3D}
\end{figure*}

\begin{figure}[ht]
\centering
	\includegraphics[width=1\textwidth]{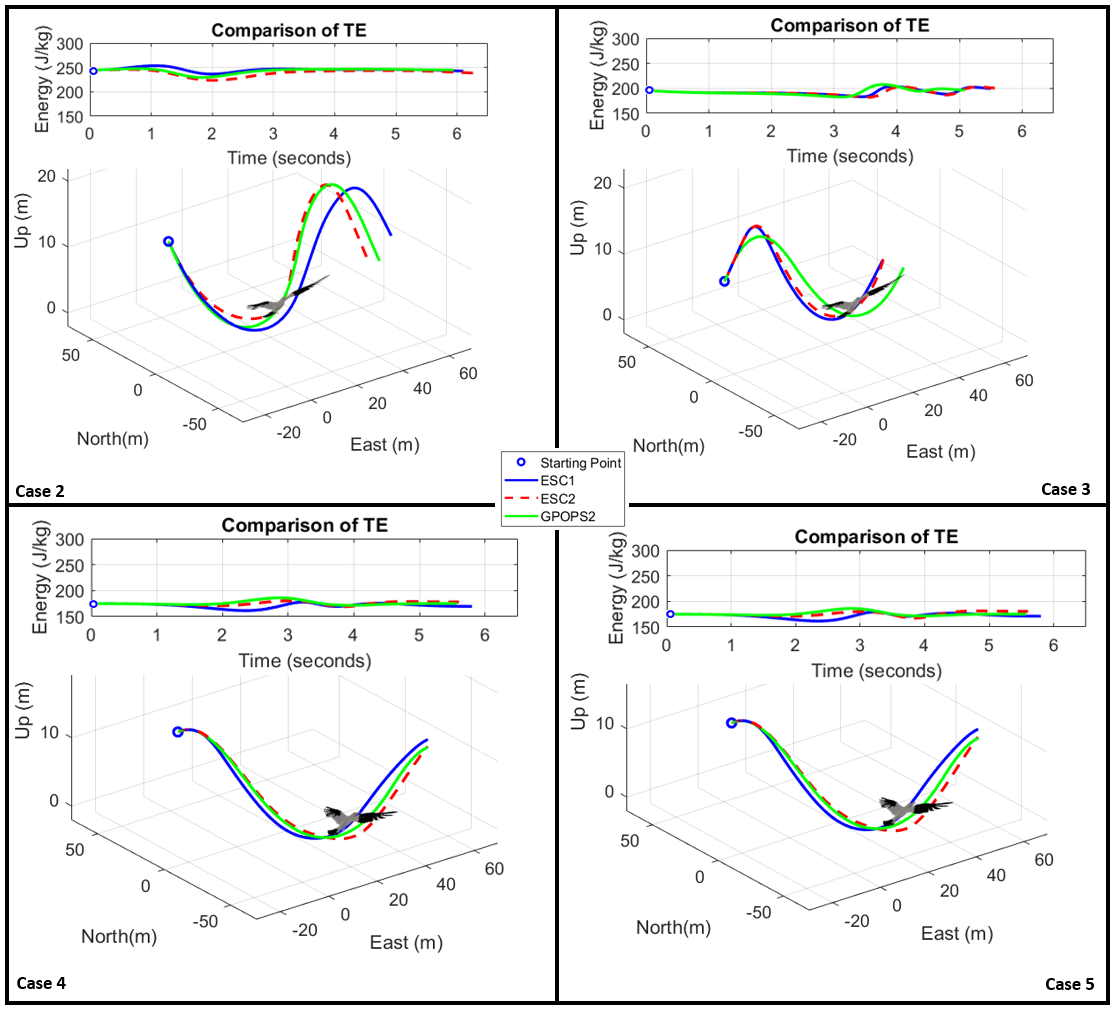}
	\caption{Case 1-5 simulation results. Depicted above is a comparison of the three dimensional trajectories and total energy of an albatross/mimicking-object  conducting dynamic soaring using extremum seeking systems (ESC1, ESC2), and the powerful optimal control solver, GPOPS2. Here, in case 2 (top-left) the wind conditions were assumed to follow a logarithmic wind model, where the objective function (performance index) is aiming at maximizing energy gain. In case 3 (top-right) the wind conditions, however, are following logistic wind model and the objective function is aiming at maximizing the  maximum total energy of the system. Case 4 (bottom-left) uses logarithmic wind  model with the objective function aiming at maximizing the total energy. Finally, in case 5 (bottom-right), logarithmic wind model was used with the objective function aiming at maximizing the total energy while accounting to disturbance/noise $\eta_{dist}$ computations inaccuracies or sensors noise. The simulations depicted in this figure provide a crucial evidence for the legitimacy of our novel hypothesis: dynamic soaring can be seen as manifestation of extremum seeking systems in nature. Furthermore, the depicted simulations suggest that dynamic soaring can be modeled, replicated, simulated, and controlled/mimicked in a stable and real-time basis similar to extremum seeking systems.}
	\label{fig:case2_5}
\end{figure}

We ran the simulation for one cycle of dynamic soaring. It is important to emphasize here that extremum seeking control systems -- as it is the case in most autonomous and automatic controllers -- require tuning of parameters, which is a challenge. So, for the choice of parameters in the simulations, we have followed the guidelines presented in \cite{ariyur2003real} and tuned them for optimal results. These parameters as well as initial conditions for each simulation are provided in the appendix \ref{sec:app_design_param}.
Next, we present the results of our simulation implementing ESC1, ESC2, and GPOPS2 in figures \ref{fig:case1_allenergies}- \ref{fig:case2_5}. We divide the results into 5 different cases based on the used wind profile, performance index, and noise/disturbance. The first four cases are simulated without the use of disturbance whereas case 5 includes disturbance/noise to demonstrate the controller's ability to counter disturbances. Each of the cases is explained below.

Case 1: Wind profile - Logistic, Performance index - Energy gain ($J_1$).
Results for this case are provided in figures \ref{fig:case1_allenergies} and \ref{fig:case1_3D}. Figure \ref{fig:case1_allenergies} represents the most important indication of the success of our implementation. That is, the absolute and most important characteristic of the dynamic soaring phenomenon is the energy neutrality. In other words, the dynamic soaring itself is achieved if the kinetic energy is traded off with the potential energy and vice versa without any loss to drag, i.e., the system is a conservative-like system, unlike typical flight dynamic systems which are dissipative systems due to drag. Figure \ref{fig:case1_allenergies} also provides a comparison between potential energy (PE), kinetic energy (KE), and total energy (TE) obtained during the dynamic soaring cycle using all of the three methods considered in the implementation. It can be observed that the total energy is near-constant in all of the three methods. Figures \ref{fig:case1_3D} shows the comparison between 3D trajectories, velocity, altitude and the control input $\phi$ obtained using all three methods. The results obtained are very close to each other and are comparable. Furthermore, extremum seeking control systems (ESC1 and ESC2) did not use any constraints or bounds given in (\ref{eqn:constraints2}) and (\ref{eqn:constraints3}), and work in real-time unlike GPOPS2. The results suggest that dynamic soaring is naturally manifested when extremum seeking is used as the modeler and controller. We find similar results and observations in the next 4 cases as well. For rest of the cases, we will provide plots comparing total energy, and 3D trajectories obtained by the use of all three methods ESC1, ESC2, and GPOPS2.

%%%%%%%%%%%%%%%%%%% sigmoid Gain
% \begin{figure*}[ht]
% \centering
% 	\includegraphics[width=\textwidth]{Energy_sigmoid_gain.jpg}
% 	\caption{Comparison of total energy during a cycle of dynamic soaring obtained using ESC1, ESC2, GPOPS2 with logistic wind profile and $J_1$ as the performance index. The results obtained from all three methods are similar. Also, the total energy is near-constant in all of the three methods which indicates the success of our implementation.}}
% 	\label{fig:energy11}
% \end{figure*}
% \begin{figure}[ht]
% \centering
% 	\includegraphics[width=0.75\textwidth]{3D_sigmoid_gain.jpg}
% 	\caption{Comparison of trajectories of albatross obtained using ESC1, ESC2, GPOPS2 with logistic wind profile and $J_1$ as the performance index. The trajectories obtained with different methods are similar to each other.}}
% 	\label{fig:3D11}
% \end{figure}

% \begin{figure}[ht]
% \centering
% 	\includegraphics[width=0.75\textwidth]{zvphi_sigmoid_gain.jpg}
% 	\caption{ Comparison of velocity, altitude and control input $\phi$ of albatross obtained using ESC1, ESC2, GPOPS2 with logistic wind profile and $J_1$ as performance index. The evolution of the parameters with time is comparable even though they are obtained from three different methods.}}
% 	\label{fig:zVPhi11}
% \end{figure}

%%%%%%%%%%%%%%%%%%%%%% Log gain

{Case 2: Wind profile - Logarithmic, Performance index - Energy gain ($J_1$).}
Results for this case are provided in figure \ref{fig:case2_5} (top-left). The total energy obtained using all these methods is constant or near constant, which shows the system is conservative (or conservative-like) in terms of energy. Similarly, the figure shows the comparison between 3D trajectories using all three methods, ESC1, ESC2, and GPOPS2. As it was in the previous case, the results are very close and comparable.
% \begin{figure*}[ht]
% \centering
% 	\includegraphics[width=\textwidth]{Energy_log_gain.jpg}
% 	\caption{Comparison of total energy during a cycle of dynamic soaring obtained using ESC1, ESC2, GPOPS2 with logarithmic wind profile and $J_1$ as the performance index.}
% 	\label{fig:energy21}
% \end{figure*}
% \begin{figure}[ht]
% \centering
% 	\includegraphics[width=0.75\textwidth]{3D_log_gain.jpg}
% 	\caption{Comparison of trajectories of albatross conducting dynamic soaring using ESC1, ESC2, GPOPS2 with logarithmic wind profile and $J_1$ as the performance index.}
% 	\label{fig:3D21}
% \end{figure}
%%%%%%%%%%%%%%%%%%%%%%%% Sigmoid TE

{Case 3: Wind profile - Logistic, Performance index - Total energy ($J_2$).}
Figure \ref{fig:case2_5} (top-right) shows the results in this case. Similar to the previous cases, the figure suggests that dynamic soaring is achieved, as the total energy is near-constant, which means the kinetic energy is traded off with the potential energy and vice versa without any loss to drag. Similarly, the figure shows the comparison between 3D trajectories using all three methods, ESC1, ESC2, and GPOPS2. As it was in the previous case, the results are very close and comparable.
% \begin{figure*}[ht]
% \centering
% 	\includegraphics[width=\textwidth]{Energy_sigmoid_TE.jpg}
% 	\caption{Comparison of total energy during a cycle of dynamic soaring obtained using ESC1, ESC2, GPOPS2 with logistic wind profile and $J_2$ as the performance index}
% 	\label{fig:energy12}
% \end{figure*}
% \begin{figure}[ht]
% \centering
% 	\includegraphics[width=0.75\textwidth]{3D_sigmoid_TE.jpg}
% 	\caption{Comparison of trajectories of an albatross conducting dynamic soaring using ESC1, ESC2, GPOPS2 with logistic wind profile and $J_2$ as the performance index.}
% 	\label{fig:3D12}
% \end{figure}

%%%%%%%%%%%%%%%%%%%% Log TE

{Case 4: Wind profile - Logarithmic, Performance index - Total Energy ($J_2$).}
Results for this case are provided in figures \ref{fig:case1_allenergies} (bottom-left). The figure shows the comparison between 3D trajectories using all three methods, ESC1, ESC2, and GPOPS2. As it was in the previous case, the results are very close and comparable.
% \begin{figure*}[ht]
% \centering
% 	\includegraphics[width=\textwidth]{Energy_log_TE.jpg}
% 	\caption{Comparison of total energy during a cycle of dynamic soaring obtained using ESC1, ESC2, GPOPS2 with logarithmic wind profile and $J_2$ as performance index.}
% 	\label{fig:energy22}
% \end{figure*}
% \begin{figure}[ht]
% \centering
% 	\includegraphics[width=0.75\textwidth]{3D_log_TE.jpg}
% 	\caption{Comparison of trajectories of albatross obtained using ESC1, ESC2, GPOPS2 with logarithmic wind profile and $J_2$ as performance index.}
% 	\label{fig:3D22}
% \end{figure}

Case 5: Wind profile - Logarithmic, Performance index - Total Energy ($J_2$) with Disturbance ($\eta_{dist}$). Extremum seeking control systems are characterized as input-output systems and can tolerate disturbance and noise at both ends (the input and the output) \cite{ESC2_2,stochasticPerterb}; they also can deal with uncertainty in the models and parameters (see applications in  \cite{ariyur2003real}). In fact, the augmented extremum seeking structure from which we obtained ESC2 has been studied and characterized for longtime in literature as being able to reject disturbances/noises. Said disturbances/noises are added to the measurement of the objective function since extremum seeking systems can operate in a model-free basis and the output is the feedback itself of the system \cite{ESC2_2,ariyur2003real} -- recall also the discussions of section \ref{ESC}. In ESC2, the input is the rolling control $\phi$ and the output is the objective function $J$. In this simulation case, we focus on disturbances that are reflected in the output $J_2$ (maximization of the total energy). These disturbances can be due to gusts/disturbances in the wind. Additionally, for one to obtain measurements of the objective function $J_2$, one needs to measure the velocity and altitude. This may be done via sensors that admit noises or through computations that admit some inaccuracies. Hence, for this simulation case, we added a uniform random disturbance/noise to the measurement of the performance-index/objective-function (the output of the system) to subsume inaccuracies arising in the system due to the computations or the sensors. The disturbance $\eta_{dist}$ ranges from $-5\%$ to $5\%$ of the measured value of the performance index. On the other hand, for classic extremum seeking systems (ESC1), while there is no particular $\eta_{dist}$ included in figure \ref{fig:ESC_Classic_custom}, we are motivated by the fact that classic extremum seeking structures are robust controllers \cite{TanHistory2010}, so, we added noise in ESC1 as well. Now, we compare the result between the classic structure (ESC1), augmented structure (ESC2) with added disturbance/noise, and GPOPS2 without adding any disturbance.
The initial condition and parameters used for the simulation of this case are the same as that of case 4 and the results for this case are provided in figure \ref{fig:case2_5} (bottom-right). Results show that extremum seeking systems are capable of rejecting the disturbance. The figure provides a comparison between total energy and 3D trajectories obtained during the dynamic soaring cycle using ESC1, ESC2 and GPOPS2 with the disturbance added to ESC1 and ESC2. The results confirms the ability of ESC1 and ESC2 to reject disturbances.

% \begin{figure*}[ht]
% \centering
% 	\includegraphics[width=\textwidth]{Energy_log_TE_with_noise.jpg}
% 	\caption{Comparison of total energy during a cycle of dynamic soaring obtained using ESC1, ESC2, GPOPS2 with logarithmic wind profile and $J_2$ as performance index with uniform random noise $\eta_{noise} = \pm5 \%$ of the value of the performance index added to the performance index. Results show that the system is able to reject }
% 	\label{fig:energy22_noise}
% \end{figure*}
% \begin{figure}[ht]
% \centering
% 	\includegraphics[width=0.75\textwidth]{3D_log_TE_with_noise.jpg}
% 	\caption{Comparison of trajectories of albatross obtained using ESC1, ESC2, GPOPS2 with logarithmic wind profile and $J_2$ as performance index with uniform random noise $\eta_{noise} = \pm5 \%$ of the value of the performance index added to the performance index.}
% 	\label{fig:3D22_noise}
% \end{figure}

Remark on the simulations: In all simulations of cases 1-5, while we used the same performance indices (objective functions), handling of ESC1 and ESC2 of them differs from GPOPS2. That is, extremum seeking controls work in real-time, so, the value of the performance index used for the optimization is instantaneous. However, GPOPS2 does not work in real-time, which allows it to find the optimized trajectory for the whole  simulation process.

\subsection{Discussion on extremum seeking implementation when compared to optimal control methods found in literature}
% Even though the results are similar in achieving energy neutrality, the application of extremum seeking} represents strong and novel benefits to the dynamic soaring problem. For example, the controller is autonomous and acts in a real-time, with no need for the constraints and bounds in (\ref{eqn:constraints1}),(\ref{eqn:constraints2}) and (\ref{eqn:constraints3}), unlike the numerical optimizers GPOPS2 and DC. Additionally, the controller is much simpler and appears to be very natural to this nature-inspired phenomenon. As a summary of the difference in control characteristics, a comparison between extremum seeking} structures and numerical optimizers such as GPOPS2 is presented in table \ref{tab:Comparison}.

Here, we discuss the characteristics of extremum seeking control systems when compared to optimal control solvers found in literature based on various criteria.

{Computational time:} The control works found in dynamic soaring literature are almost exclusively configured as optimal control problems in part or in full -- recall figure \ref{fig:traj_plan}. Said optimal control solutions are based on a given initial condition and predetermined wind profile models. They need a certain amount of time to obtain the optimal solution. Thus, they can not be implemented in real-time. Extremum seeking controllers, on the other hand, work in real-time. By this we mean, ESC automatically generates the control trajectories as it progresses (marches) with time, based on the measurement of the objective function taken as feedback.

{Dependency on the model:} Extremum seeking systems are model-free adaptive control techniques, whereas optimal control solvers used to solve dynamic soaring in literature are model-dependent. For instance, optimal control solvers (refer section \ref{problem_setup}) make use of the mathematical expression of the performance index (objective function) to perform further mathematical operations/computations, such as but not limited to, the computation of the gradient, and in that case, the closed form of the mathematical expression is required. On the other hand, for extremum seeking systems, the gradient of the objective function does not need to be calculated, as the extremum seeking structure itself drives the system to its optimized objective via perturbing the measured/sensed objective function, which causes automatic realization of the gradient. This exempts the need for the closed form of the objective function and only the measured value of the performance index is sufficient. Furthermore, taking into the consideration that wind shear models themselves are an approximation of actual wind conditions, and the entire dynamic soaring model can be less than proper, the overall accuracy of optimal control solvers can be compromised at times. But for extremum seeking systems, measurement of the objective function (even with some noise as shown in case 5 in section \ref{simulationSection}) is enough for the controller to work properly. Hence, extremum seeking structures are less dependent on the model of the system when compared to optimal control solvers in literature.

{Parameters and constraints:} Optimal control solvers need to make use of many constraints and bounds to obtain the optimal solutions, however, these are not required for extremum seeking systems. In other words, for basic extremum seeking system implementation, the bounds and constraints in (\ref{eqn:constraints2})-(\ref{eqn:constraints3}) are not a prerequisite for the implementation.
It is worth noting that some of these constraints used in optimal control solvers are safety constraints that can guide the aircraft to avoid dangerous situations. Indeed, safety constraints are essential for the safe flight of an unmanned system or a mimicking-object performing dynamic soaring. While the work of this paper only aims at providing proof of concept that dynamic soaring itself emerges when we use extremum seeking as the controller without having the need of any constraints,
% However, the point we are trying to prove via this proof of concept work is that dynamic soaring itself emerges when we use extremum seeking as the controller without having the need of the constraints to obtain the path. 
having constraints as safety bounds should not contradict the use of extremum seeking controllers. As a matter of fact, the extremum seeking literature itself has methods that operate with constraints \cite{liao2019constrainedESC}. So, in future research, these constraints can be added if needed from the safety point of view.
As for the parameters, proper tuning is still a challenge in extremum seeking as it is with other techniques found in automatic control literature, e.g., in proportional-integral-derivative (PID) control designs. However, once properly tuned, the extremum seeking system can work in real-time. Optimal control solvers on the other hand also have multiple parameters. For instance, one has to assume the appropriate time-steps size, whether the time steps and the grids of the parameter space should be uniform or adaptive, and guess values that need to be imposed on the system, etc. Most importantly, optimal control solvers do not work in real-time with different grades of complexity.

{Stability of the  implemented extremum seeking} structures:
As discussed in the case 5 of the simulations earlier, extremum seeking control systems are characterized as input-output systems and can tolerate disturbance and noise at both ends (the input and the output) \cite{ESC2_2,stochasticPerterb}. The stability of the classic extremum seeking structure as in figure \ref{fig:ESC_Classic} has been rigorously studied and characterized in \cite{ariyur2003real,KRSTICMain}. On the other hand, the stability of the augmented extremum seeking structure as in figure \ref{fig:ESCPlant} has been studied in \cite{ariyur2003real,ESC2_2}. Both of said systems/structures -- which we applied to dynamic soaring -- are stable and capable of rejecting disturbances. The disturbance rejection by the system is an important indication of stability. As seen in case 5 in section \ref{simulationSection}, both the classic structure and the augmented structures of extremum seeking (ESC1 and ESC2) were able to reject disturbances added to the measurement of the objective function, and provide dynamic soaring trajectory with near-constant total energy. This is an important indication of the success of our control design. In addition, the augmented extremum seeking structure (ESC2) provides a more customized/particular method to prove rigorously the stability of the structure when customized to a particular problem. For that, one needs to properly choose Block 1 and Block 2 in the design of the augmented extremum seeking system (as shown in figure \ref{fig:ESCPlant}). The design of the blocks depends on the input and output optional dynamics $F_i$, $F_0$ of the extended dynamical system and $\Gamma_\phi$, $\Gamma_J$, the nature/guess of the optimal parameter $\phi^*$, and the optimal objective value $J^*$. For the dynamic soaring problem in ESC2 (see figure \ref{fig:ESCPlant_custom}), we take $F_i =1,F_0=1$, and Block 1 and Block 2 are designed such that assumptions C1-C5 (provided in the appendix \ref{sec:app_algorithm}) are satisfied. This leads to the stability guarantee for the system. The reader is directed to the appendix \ref{sec:app_algorithm} for more on the assumptions, design algorithm, rigorous stability analysis, and our stability theorem.
% choose the design elements of the blocks and functions in the extremum seeking} structure. The choice of design elements can be done by following certain assumptions and the design algorithm presented in \cite{ESC2_2}. We have provided the assumptions, design algorithm, and stability analysis in the appendix \ref{sec:app_algorithm}. 

The above discussion on extremum seeking performance compared to optimal control methods found in literature is summarized in table \ref{tab:Comparison}.
 \begin{table}[ht]
\centering
\resizebox{\textwidth}{!}{
    \begin{tabular}{|>{\arraybackslash}p{5cm}|>{\arraybackslash}p{5cm}|>{\arraybackslash}p{5cm}|}
		\hline 
		 \textbf{Performance criterion} & \textbf{Extremum seeking structures} & \textbf{Numerical optimizers such as GPOPS2}\\
		 \hline
		Functionality in Real-time. & Work in real-time. & Do not work in real-time.\\
		\hline
		Nature of the objective function. & Measured value of the objective function is enough and does not necessitate the mathematical expression a priori. & Need mathematical expression of the objective function for further computation.\\
		\hline
		Design parameters, constraints, and bounds. & Some design parameters need to be tuned and do not require constraints and bounds. & Heavy amount of parameters, bounds, and constraints are needed.\\
		\hline
		Stability. & Stable. & Not necessarily Stable \cite{mir2021stability}.\\
		\hline
	\end{tabular}
	}
	\caption{Comparison between the proposed extremum seeking} structures against other numerical optimizers on various criteria.
	\label{tab:Comparison}
\end{table}

% \subsection{Determining the stability of the ESC1 structure}\label{stabilityDiscussion}
% The design of the ESC1 structure is done in a such way that the assumptions 2.1 -2.3 and the conditions required by Theorem 2.1 in \cite{ESC2_2} are satisfied. The compensator $C_0=(s \sin 2.5 + \cos 2.5)/((s+0.4)(s+2))$ is asymptotically stable, the block 1 and block 2 $ (s^2+1.5^2)/((s+0.4)(s+2))$ are proper, $\Gamma_J$ and $\Gamma_\phi$ are strictly proper. 
% Before doing stability analysis, we define the following terms. 
% 
    % $C_0(s)$ and $1/(1+L(s))$ are asymptotically stable, where $L(s)=af''/4 \cdot H_i(s) (e^{j\phi_0} F_i(j\omega) H_0(s+j\omega))+ e^{-j\phi_0}F_i(-j\omega)H_0(s-jw)$, and $H_i=C_i(s) \Gamma_\phi (s), F_i(s)$,  $H_0(s)=C_0(s)/ \Gamma_J(s)F_0(s)$.

    % Let the smallest in absolute value among the real parts of all of the poles of $H_{osp}(s)$ be denoted by $a$. Let the largest among the moduli of all of the poles of $F_i(s)$ and $H_{obp}$ be denoted by $b$, where $H_0(s)=H_{osp}H_{obp}=kC_o(s)/\Gamma_J(s)$. The ratio $M=a/b$ is sufficiently large.

\section{Conclusive remarks and future works}\label{conclusion}
This paper provides a novel hypothesis on how soaring birds -- such as the albatross -- perform the energy-efficient dynamic soaring maneuver by utilizing energy from the wind -- similar to the process depicted in figure \ref{fig:dynamicsoaring}; in summary, the paper investigated whether the albatross is performing dynamic soaring in a similar manner to extremum seeking systems as philosophized in figure \ref{fig:philosophy}. The paper has provided a proof of concept for said hypotheses by successfully characterizing, modeling, and controlling dynamic soaring as an extremum seeking system with successful simulations and comparisons with the literature. Moreover, the paper provides a novel contribution framework and original insights for future research in two fundamental ways. The first fundamental is regarding the dynamic soaring problem itself and the enabling of bio-mimicry of the albatross, and other soaring birds, by unmanned systems. The results of this paper, while having limitations, provide a less-constrained, real-time, autonomous, stable, and simple control solution for the dynamic soaring problem. This paper simply shows that dynamic soaring can be interpreted as a manifestation of extremum seeking systems in nature and can be captured, replicated, and performed by extremum seeking systems; this opens the door for dynamic soaring to benefit from the rich literature of extremum seeking control systems in both the theoretical and implementation fronts. The second fundamental is regarding the topic of extremum seeking in nature. Arguably, the process and mechanisms many birds, animals, and insects are conducting to reach optimally (through dynamic optimization) their objective, can be seen as an extremum seeking system mechanism. Organisms in nature cannot be conducting non-real-time, unstable, or very complex procedures due to the lack of resources and computational ability. We believe this paper contributes to the legitimacy of the above-mentioned observation that many of organisms in nature might be conducting extremum seeking processes. 
In the future, we are aiming to advance our research in two major ways. First, we are aiming at implementing the dynamic soaring problem in a more complex situation. For instance, we aim at implementing multiple extremum seeking loops incorporating multi-variable control inputs to allow pitching and rolling simultaneously. Also, we aim at investigating the application of extremum seeking systems when the system faces poor wind conditions and need to complement the energy from the wind by thrust. Moreover, additional research needs to be done on sensing and how one would take measurements to enable extremum seeking implementation. Experimental validation of this paper's hypothesis will be a major route for future research as well. Second, we aim at investigating other optimized phenomena in systems biology and nature, especially those that are associated with oscillatory/periodic-like behaviors as this may indicate possible parallelism between such phenomena and extremum seeking processes.
We believe this paper will encourage us and others to pursue more investigations of natural phenomena that can be characterized by extremum seeking systems.        
\section{Acknowledgement}
We would like to thank the member of the editorial board of Bioinspiration \& Biomemetics who, in addition to the invited reviewers, has reviewed this paper. This member has provided very helpful comments on the paper, especially for its readability and how it can be improved for interdisciplinary readership to target different communities of science and engineering who can have a strong interest in the topic of this paper. The authors are indebted to this member; the provided comments by this member were transformational in how this paper appears at the moment.

\section{Appendix}
\subsection{Design and stability analysis of the augmented extremum seeking structure}\label{sec:app_algorithm}
 
 We will first expand on the reasoning behind the selection of the nature/guess of $\phi^*(t)$ and $J^*(t)$ as periodic or sinusoidal-like for the augmented extremum seeking structure (ESC2), then we will describe the process to design the structure itself, and finally analyze its stability.

In the controllability study published recently on dynamic soaring \cite{mir2018controllability}, it was shown that the dynamic soaring control system in \eqref{eqn:probForm} can be formulated as a control-affine system of the form \eqref{eqn:controlAffine}: 
\begin{equation}\label{eqn:controlAffine}
    \dot{\bm{x}}=\bm{b_0}(\bm{x})+\sum_{i=1}^m u_i \bm{b_i}(\bm{x}),
\end{equation}
where $\bm{x}\in\mathbb{R}^n$ is the state space vector, $\bm{b}_0$ is the drift (uncontrolled) vector field of the system, $u_i$ are the control inputs (pitch rate and roll rate) , and $\bm{b}_i(\bm{x})$ are the control vector fields. Given that we are dealing with a single-input single-output system, the only control input variation we can apply will be through the rate of $\phi$ (roll rate). The controllability of the system \eqref{eqn:controlAffine} require control inputs that satisfy the following condition: if $U \in \mathbb{R}^m$ is the set of admissible controls then for every $u_i$, a given control choice in the admissible set, we have $u_i:I\rightarrow U$ for $I\subset \mathbb{R}$ locally integrable and also $\bm{0}\in \text{int}(\text{conv}(U))$ -- see chapter 7, \cite{bullo2004geometric}, where $\text{int}(\text{conv}(U))$ refers to the interior of the convex hull of $U$. This condition needs to be satisfied in order for the so called ``Lie bracket proper excitation/variation" to take place. Now, as shown in \cite{sussmann1993lie, DURR2013}, sinusoidal inputs are typical candidates to cause the control inputs variation/excitation which cause steering/motion-planning of the system and provide controllability. Hence, the optimal input and output can be reasonably guessed as periodic and sinusoidal-like. This is supported by many observations made on the periodicity of dynamic soaring maneuvers, theoretically as in \cite{mir2018optimal,gao2017dubins,mir2019soaring,MITbousquet2017optimal} and experimentally as in \cite{sachs2013experimental,yonehara2016flight}.

Having now made a guess on the nature of the optimal input $\phi^*(t)$ and the optimal output $J^*(t)$ and computed their Laplace transform $\Gamma_{\phi}$ and $\Gamma_J$ as in (\ref{eqn:laplace_transform}), next we present an algorithm taken from \cite{ESC2_2} as a guideline for designing the compensators $C_i$ and $C_0$. Before describing the design process, we would like to remind the reader of some technical terms associated with the design process. The poles and zeros of a transfer function are the roots of its denominator and numerator, respectively. Similarly, a transfer function is called ``a proper transfer function" if the degree of its numerator does not exceed the degree of its denominator. Now, we describe the design process for the augmented extremum seeking structure (ESC2) in figure \ref{fig:ESCPlant_custom}. For step 1, we need to select sufficiently large frequency $\omega$ but should avoid choosing it equals to any frequency persistent in the disturbance $\eta_{dist}$ as this may lead to large steady state tracking error. Generally, if $\pm j\omega$ is equal to a zero of the input dynamics $F_i(s)$, then the dynamics of the system will not be perturbed properly. However, in our case, we have $F_i(s)=1$ (recall the discussion on this from section \ref{ESC_DS_structure}), so, this has been avoided.  Similarly, for step 2, the choice of the amplitude $a$ should be small enough to have a small magnitude of the input perturbation $a|F_i(j\omega)|$, but large enough to have a measurable variation in the output. Again, we only focus on the amplitude $a$ since $F_i$=1. Now, in step 3, for the design of $C_0(s)$, three important criteria should be considered: (i) $C_0(s)$ should have zeros of $\Gamma_J(s)$ that are not asymptotically stable as its zeros; (ii) the compensator $C_0(s)$ itself should be asymptotically stable; and (iii) the ratio $\frac{C_0(s)}{\Gamma_J(s)}$ should be proper.
Finally in step 4, for the design of $C_i(s)$, three important criteria should be adopted as well: (i) $C_i(s)$ should not include any poles of the $\Gamma_{\phi}(s)$ that are not asymptotically stable as its zeros; (ii) the product $C_i(s)\Gamma_{\phi}(s)$ should be proper; and (iii) $1/(1+L(s))$ should be asymptotically stable, where $L(s)$ is defined later in the assumption C5. These guidelines are summarized in Algorithm \ref{alg:alg1}.

\begin{algorithm}
		\caption{Algorithm for compensator design} 
% 		\color{blue}
		\label{alg:alg1}
			\begin{algorithmic}[1]
			\State Select the perturbation frequency $\omega$ sufficiently large and not equal to any frequency in the disturbance $\eta_{dist}$ and with $\pm j\omega$ not equal to any imaginary axis zero of $F_i(s)$.

			\State Set the perturbation amplitude ``$a$" to obtain small steady state output error $\tilde{J}(t)=J(t)-J^*(t)$.
			
			\State Design $C_0(s)$ asymptotically stable with zeros of $\Gamma_J(s)$ that are not asymptotically stable as its zeros, and such that $\frac{C_0(s)}{\Gamma_J(s)}$ is proper. 
			
			\State Design $C_i(s)$ by any linear single-input single-output design technique such that it does not include poles of $\Gamma_{\phi}(s)$ that are not asymptotically stable as its zeros, $C_i(s)\Gamma_{\phi}(s)$ is proper, and $1/(1+L(s))$ is asymptotically stable.
			
		\end{algorithmic} 
\end{algorithm}

Now, we rigorously prove that the design we have chosen for the augmented extremum seeking system structure for dynamic soaring (ESC2) in figure \ref{fig:ESCPlant_custom} is stable. However, before providing our stability theorem, we introduce the following assumptions and conditions C1-C5 as in \cite{ESC2_2}:
\begin{enumerate}[label= C\arabic*.]
    \item 
    $F_i(s)$ and $F_0(s)$ are asymptotically stable and proper.
    \item
    $\Gamma_J(s)$ and $\Gamma_\phi(s)$ are strictly proper rational functions and the poles of $\Gamma_\phi(s)$ that are not asymptotically stable are not zeros of $C_i(s)$.
    \item
    Zeros of $\Gamma_J(s)$ that are not asymptotically stable are also zeros of $C_0 (s)$.
    \item
    $C_0(s)/\Gamma_J(s)$ and $C_i(s) \Gamma_\phi (s)$ are proper.
    \item
    $C_0(s)$ and $1/(1+L(s))$ are asymptotically stable, where $L(s)=af''/4 \cdot H_i(s) [Re(e^{j\phi_0} F_i(j\omega) H_0(s+j\omega))]$ and $H_i=C_i(s) \Gamma_\phi (s) F_i(s)$,  $H_0(s)=C_0(s)/ \Gamma_J(s) \cdot F_0(s)$, $f''$ is constant.
\end{enumerate}
Now, we present our theorem for stability as

\textbf{Theorem 1:} The output error $\Tilde{J}=J-J^*$ in the system ESC2, depicted in figure \ref{fig:ESCPlant_custom}, achieves local exponential convergence to an $O(a^2+1/\omega^2)$ neighborhood of origin.

\textit{Proof}: 
Since $F_i(s)=1$ and $F_0(s)=1$ are asymptotically stable and proper, the condition C1 is satisfied. Similarly, $\Gamma_J$ and $\Gamma_\phi$ as defined in (\ref{eqn:gammas}) are strictly proper rational functions. So, the condition C2 is also satisfied. In addition, with choosing $C_i(s)=1$ and $C_0(s)=(c_3 \cos c_0 + s \sin c_0)/((s+c_4)(s+c_5)(s+c_6))$, we get $C_i$, $C_i(s)\Gamma_\phi $, and $ C_0(s)/\Gamma_J $ that satisfy the conditions C2-C5. Now, from Theorem 2.1 in \cite{ESC2_2}, satisfaction of conditions C1-C5 completes the proof.

{Remark}: The structure in figure \ref{fig:ESCPlant_custom} is a member of a family of reduced structures from the generalized structure in figure \ref{fig:ESCPlant}. Hence, the proof presented in Appendix B of \cite{ESC2_2} applies to the structure in  figure \ref{fig:ESCPlant_custom} as well.

\subsection{Design parameters for various cases} \label{sec:app_design_param}
Here, we provide the parameters we used for the design and implementation of ESC1 and ESC2 structures as well as the initial conditions used in each of the simulation cases ( section \ref{simulationSection}). Here, each case corresponds to the use of different wind shear model and performance indices. Case 1 refers to the use of logistic wind model and maximum energy gain as performance index. Similarly, case 2 refers to the use of logarithmic wind model and maximum energy gain as performance index. Case 3 is corresponding to the use of logistic wind model and maximum total energy as performance index. Similarly, case 4 refers to the use of logarithmic wind model with maximum total energy as the performance index. Finally, case 5 refers to the use of logarithmic wind model with maximum total energy as performance index and with disturbance $\eta_{dist}$. Table \ref{tab:ESC1_params} provides the parameters for ESC1 structure and table \ref{tab:ESC2_params} provides the parameters corresponding to ESC2 structure. The initial state of the simulation $\bm{x_0} = [x,y,z,V,\gamma,\psi]$ is provided in table \ref{tab:Initial_condition} for all of the cases.

\begin{table}[h]
\centering
\resizebox{0.6\textwidth}{!}{
    \begin{tabular}{|>{\arraybackslash}p{1.5cm}|>{\arraybackslash}p{0.5cm}|>{\arraybackslash}p{0.5cm}|>{\arraybackslash}p{0.5cm}|>{\arraybackslash}p{1.5cm}|>{\arraybackslash}p{0.5cm}|>{\arraybackslash}p{0.5cm}|}
\hline
Case & $\omega_1$ &$a_1$ & $b_1$& $\phi_{phase_1}$ &  $k_1$ &  $h$ \\
\hline
 1 & 1.2&	0.5&	0.2	&0.1&	0.1	&0.4 \\
\hline
2&	5.0&	0.2&	0.8&	-1.9	& 1.8 &	0.4\\
\hline
3 &1.9&	0.6&	1.6&	0.6&	0.1&	4.8\\
\hline
4  &0.5&	0.6&	0.3&	-0.2&	0.7&	1.0\\
\hline
5  &0.5&	0.6&	0.3&	-0.2&	0.7&	1.0\\
\hline

\end{tabular}
}
\caption{Parameters used in ESC1 structure in various simulation cases.}
\label{tab:ESC1_params}
\end{table}

% For augmented structure
\begin{table}[h]
\centering
\resizebox{\textwidth}{!}{
    \begin{tabular}{|>{\arraybackslash}p{1.5cm}|>{\arraybackslash}p{0.5cm}|>{\arraybackslash}p{0.5cm}|>{\arraybackslash}p{0.5cm}|>{\arraybackslash}p{1.5cm}|>{\arraybackslash}p{0.5cm}|>{\arraybackslash}p{0.5cm}|>{\arraybackslash}p{0.5cm}|>{\arraybackslash}p{0.5cm}|>{\arraybackslash}p{0.5cm}|>{\arraybackslash}p{0.5cm}|>{\arraybackslash}p{0.5cm}|}
\hline
Case & $\omega_2$ &$a_2$ & $b_2$& $\phi_{phase_2}$ &  $k_2$ & $c_1$  & $c_2$ &  $c_3$ & $c_4$ & $c_5$ & $c_6$   \\
\hline
 1 & 1 & 0.4 & 1.8  & -0.8   & 1.5& 8.2& 1.8 & 1.5 &0.1 & 8.8 & 8.1  \\
\hline
2&	0.8&	0.8&	1.5&	-2.2&	1.3&	2.3&	9.3&	1.0&	0.4&	3.5&	3.0\\
\hline
3 &	1.7&	0.2&	0.5&	0.5&	1.4&	0.7&	3.1&	1.5&	1.8&	9.7&	2.2\\
\hline
4  &	2.4&	0.7&	1.7&	-2.9&	1.3&	3.8&	6.3&	1.1&	3.1&	9.8&	9.6\\
\hline
5  &	2.4&	0.7&	1.7&	-2.9&	1.3&	3.8&	6.3&	1.1&	3.1&	9.8&	9.6\\
\hline

\end{tabular}
}
\caption{Parameters used in ESC2 structure in various simulation cases.}
\label{tab:ESC2_params}
\end{table}

% For initial coniditio
\begin{table}[h]
\centering
\resizebox{0.6\textwidth}{!}{
    \begin{tabular}{|>{\arraybackslash}p{1.5cm}|>{\arraybackslash}p{1cm}>{\arraybackslash}p{1cm}>{\arraybackslash}p{1cm}>{\arraybackslash}p{1.5cm}>{\arraybackslash}p{1cm}>{\arraybackslash}p{1cm}|}
\hline
Case &  & &$\bm{x}_0$ & &  &   \\
\hline
 1 & -16.0&	15.0&	10.0&	14.0&	-0.7&	-0.1 \\
\hline
2&	-16.0&	15.0&	15.0&	14.0&	-1.0& -0.6\\
\hline
3 &-16.0&	15.0&	10.0&	14.0&	0.6&	0.4\\
\hline
4  &-16.0&	15.0&	15.0&	7.5&	0.4&	0.3\\
\hline
5  &-16.0&	15.0&	15.0&	7.5&	0.4&	0.3\\
\hline
\end{tabular}
}

\caption{Initial condition for various simulation cases.}
\label{tab:Initial_condition}
\end{table}

\section*{References}
% \begin{verbatim}
% \documentclass[12pt]{iopart}
% \usepackage{CJK}
% .
% .
% .
% \begin{document}
% \begin{CJK*}{GBK}{ }

% \title[]{Title of article}
% \author{Author Name (CJK characters)}
% \address{Department, University, City, Country}
% .
% .
% .
% \end{CJK*}
% \end{verbatim}
% \bibliographystyle{IEEEtran}
% % \bibliographystyle{agsm} % Harvard style as given in guidelines
\bibliography{references}
\bibliographystyle{IEEEtran}
% \bibliography{bibliography}

\end{document}